\documentclass[9pt,twocolumn,twoside]{IEEEtran}
\usepackage{bbm}
\usepackage{amsmath}
\usepackage{amsfonts}
\usepackage{mathrsfs}
\usepackage{amssymb}
\usepackage{multirow}
\usepackage{latexsym}
\usepackage{psfrag}
\usepackage{graphicx}
\usepackage{cite}
\usepackage{color}
\usepackage{times,balance}
\usepackage{algorithm,multirow,xcolor}
\usepackage{algorithmicx}
\usepackage{algpseudocode}
\usepackage{longtable}
\usepackage{booktabs}
\usepackage{colortbl}
\usepackage{bm}
\usepackage{booktabs}
\usepackage{amssymb}
\usepackage{pifont}
\usepackage{wasysym}
\usepackage{threeparttable}
\usepackage{subfigure}
\usepackage{caption}

\newcommand{\cmark}{\ding{51}}%
\newcommand{\xmark}{\ding{55}}%
\definecolor{mygray}{gray}{.8}

\newtheorem{theorem}{Theorem}

\newtheorem{assumption}{Assumption}

\newtheorem{lemma}{Lemma}
\newtheorem{remark}{Remark}

\allowdisplaybreaks[4]

\begin{document}

\title{Quantized Distributed Gradient Tracking Algorithm with Linear Convergence in Directed Networks}


%
%
\author{Yongyang Xiong, 
Ligang Wu, \IEEEmembership{Fellow, IEEE},  
Keyou You, \IEEEmembership{Senior Member, IEEE},
and Lihua Xie, \IEEEmembership{Fellow, IEEE}
\thanks{This work was supported by the Technology and Innovation Major Project of the Ministry of Science and Technology of China under Grant 2020AAA0108400 and Grant 2020AAA0108403,    
the National Natural Science Foundation of China (62033005, 62203254), 
and the Natural Science Foundation of Heilongjiang Province (ZD2021F001). 
\emph{(Corresponding author: Ligang Wu.)}
}
\thanks{Y. Xiong and K. You are with the Department of Automation, 
	Beijing National Research Center for Information Science and Technology, 
	 Tsinghua University, Beijing 100084, P. R. China. E-mail: \texttt{xiongyy@tsinghua.edu.cn}; \texttt{youky@tsinghua.edu.cn}.}
\thanks{L. Wu is with the
Department of Control Science and Engineering, Harbin Institute of Technology, Harbin 150001, P. R. China.
E-mail: \texttt{ligangwu@hit.edu.cn}.} 
\thanks{L. Xie is with the
School of Electrical and Electronic Engineering,
Nanyang Technological University, 50 Nanyang Avenue, Singapore 639798. 
E-mail: \texttt{elhxie@ntu.edu.sg}.}}

\maketitle

\begin{abstract}

Communication efficiency is a major bottleneck in the applications of distributed networks. 
To address the problem, the problem of quantized distributed optimization has attracted a lot of attention. However, 
most of 
the existing quantized distributed optimization algorithms can only converge sublinearly. 
To achieve linear convergence, 
this paper proposes a novel quantized distributed gradient tracking algorithm (Q-DGT) to minimize a finite sum of
local objective functions over directed networks. 
Moreover, we explicitly derive 
lower bounds for the number of quantization levels, and prove that Q-DGT can 
converge linearly even when 
the exchanged variables are respectively 
 quantized with 3 quantization levels.  
Numerical results also confirm the efficiency of the proposed algorithm.

\end{abstract}


\vspace{-1ex}
\section{Introduction}

\IEEEPARstart{R}{ecent} years have witnessed  tremendous progress in distributed optimization due to its
wide applications in formation control \cite{YouK2011Network},
distributed resource allocation \cite{ZhangJ2020Distributed},
online optimization \cite{LiX2020Distributed}, 
localization systems \cite{Xie2021Survey}, 
game theory \cite{YeM2019A},
to name a few.
They require a group of networked nodes to cooperatively optimize the sum of their local cost functions via local communications. 
A comprehensive review of 
this topic can be found in 
\cite{NedicA2018Network,NedicA2020Distributed}. 


Although distributed algorithms are capable of solving complex tasks in a collaborative manner,  
limited communication capacity is a major bottleneck in distributed networks,   
especially for large-scale distributed machine learning. How to design communication-efficient distributed algorithms has 
attracted an increasing attention \cite{Magnusson2018Convergence,DoanT2020Convergence,Doostmohammadian2022Fast}.  
%
For instance, the encoding-decoding scheme in \cite{LiT2010Distributed} has been designed to distributedly solve linear equations \cite{LeiJ2020Distributed}, distributed optimization problems \cite{YiP2014Quantized,LiH2020Computation}. 
To further reduce the size of data transmission, the recent work \cite{Zhang2019Distributed} showed that the sign of relative state between neighbors is sufficient for achieving convergence. 
As errors are inevitable for a finite-precision quantizer, 
the QDGD algorithm proposed in \cite{Reisizadeh2019An} achieves vanishing consensus error even in the presence of non-vanishing noise by modifying the contribution of the received quantized information for each node. 
%
%
%
By incorporating quantization scheme into the push-sum algorithm\cite{NedicA2009distributed}, the authors of  \cite{Taheri2020Quantized} proposed distributed algorithms over directed networks for both convex and non-convex functions.
%
Since the aforementioned works are derived by the distributed gradient descent (DGD) \cite{NedicA2009distributed}, they can only achieve sublinear convergence even for the strongly convex functions.  


How to accelerate the convergence speed is fundamentally important to reduce communication cost. 
Recently, a few significant efforts have been devoted to designing quantized distributed algorithms with linear convergence. 
For instance, 
ref. \cite{KovalevD2020A,LiuX2021Linear} proposed DQOA and LEAD, respectively, under the assumption that the randomized quantizer $Q(\cdot)$ is an unbiased and $\delta$-contracted operator, i.e. $\mathbb{E}[Q(x)]=x$ and $\mathbb{E}\Vert Q(x)-x\Vert^{2}\leq \delta \Vert x\Vert^{2}$ for all $x\in\mathbb{R}^{m}$. 
Clearly,  
this assumption 
excludes some important quantizers, e.g., the binary quantizer. 
 Ref. \cite{Lee2018Finite} proposed Q-NEXT by dynamically adjusting the center of the quantization interval. 
Ref. \cite{MagnussonS2019On} established a trade-off between the convergence speed and the communication cost 
per iteration so that linear convergence can be guaranteed. 
Although the aforementioned quantized algorithms \cite{KovalevD2020A,Lee2018Finite,LiuX2021Linear,MagnussonS2019On} converge linearly, they are designated only for undirected networks.
Note that extending distributed algorithms from undirected networks to directed networks is non-trivial \cite{XieP2018Distributed,NedicA2015Distributed,PuS2018PushJou,XinR2018A,Nedic2017Achieving}. 
In fact, if the directed network is unbalanced,  
i.e., 
	there exists at least a node that the sum of the weights of its outgoing nodes is not equal to that of its incoming nodes (see e.g., \cite{GharesifardB2012Distributed,Rikos2014Distributed}), 
the DGD finally minimizes a weighted average of local functions. 
Hence, 
an additional variable is usually exchanged between nodes to eliminate the effects of the unbalancedness\cite{NedicA2015Distributed,Nedic2017Achieving}. 
To resolve the unbalancedness issue, 
the push-pull/$\mathcal{AB}$ algorithm \cite{PuS2018PushJou,XinR2018A} and its variant \cite{ZhangJ2020Distributed}
leverage row-stochastic matrix and column-stochastic matrix simultaneously and achieve exact linear convergence for strongly convex and smooth functions. 
In sharp contrast to the subgradient-based quantized algorithms in \cite{YiP2014Quantized,LiH2020Computation}, \cite{Taheri2020Quantized} that only the decision variable needs to be quantized,  
the quantizer cannot be 
directly incorporated into 
push-pull/$\mathcal{AB}$  as it 
will result in an accumulation of quantization errors\cite{PuS2020AR}, 
thereby the convergence cannot be guaranteed.

\begin{table}[t]
\caption{\textsc{Distributed Optimization Algorithms}}
\vspace{-1em}
\begin{center}
\begin{tabular}{lcccc}
\toprule[1pt]
References ~\hspace*{-2em}& digraphs ~\hspace*{-2em}&~  linear convergence ~\hspace*{-2em}&~ 1-bit communication\\
\midrule
\cite{LiH2020Computation,Reisizadeh2019An,DoanT2020Convergence} ~\hspace*{-2em}& \xmark  ~\hspace*{-2.5em}&~ \xmark ~\hspace*{-2.5em}&~ \xmark\\
\cite{YiP2014Quantized,Zhang2019Distributed,Magnusson2018Convergence} ~\hspace*{-2em}& \xmark  ~\hspace*{-2.5em}&~ \xmark ~\hspace*{-2.5em}&~ \cmark\\
\cite{Lee2018Finite,KovalevD2020A,LiuX2021Linear,MagnussonS2019On,QuG2018Harnssing,XuJ2018Convergence}~\hspace*{-2em}& \xmark  ~\hspace*{-2.5em}&~ \cmark ~\hspace*{-2.5em}&~ \xmark\\
\cite{Taheri2020Quantized}, \cite{XieP2018Distributed,NedicA2015Distributed} ~\hspace*{-2em}& \cmark  ~\hspace*{-2.5em}&~ \xmark ~\hspace*{-2.5em}&~ \xmark\\
\cite{ZhangJ2020Distributed}, \cite{PuS2018PushJou,XinR2018A,Nedic2017Achieving} ~\hspace*{-2em}& \cmark ~\hspace*{-2.5em}&~ \cmark ~\hspace*{-2.5em}&~ \xmark\\
~Our work ~\hspace*{-2em}& \cmark  ~\hspace*{-2.5em}&~ \cmark ~\hspace*{-2.5em}&~ \cmark\\
\toprule[1pt]
\end{tabular}
\end{center}
\label{Table1}
\end{table}

\vspace{0 mm}

%
A question naturally arises: 
Whether it is possible to 
develop a quantized distributed algorithm 
over directed networks that converges linearly even for 
one-bit communication?
In this paper, we give a positive answer.
A comparison of our work with the state-of-the-art 
works 
is provided in Table \ref{Table1}. 
The main contributions of this work are summarized as follows: 
\begin{enumerate}
\item[1)]  We propose a novel quantized distributed algorithm Q-DGT over directed networks. The Q-DGT is remarkably robust to quantization errors, and achieves linear convergence. 
%


\item[2)] We explicitly provide the 
lower bounds of the quantization levels  
to resolve the saturation issue for the finite-level quantizers, which even supports the extreme 3-level quantization.

\end{enumerate}

The remainder of this paper is organized as follows. We formulate the problem in Section II. The Q-DGT is
provided in Section III. Section IV includes the convergence analysis. Simulation results are presented in Section V.
In Section VI, we conclude this paper. 

{\bf Notation.}
We use $x_{i}$ to denote the $i$-th element of vector $x$; 
$\bm{1}_{n}$($\bm{0}_{n}$) denotes a column vector with its all elements equaling to one(zero).
The notation $f=\mathcal{O}(h)$ means there exists a positive constant $\upsilon<\infty$ such that $f\leq \upsilon h$. 
$\nabla F(x(k)) \triangleq (\nabla f_{1}^{\text{T}}(x_{1}(k)),...,\nabla f_{n}^{\text{T}}(x_{n}(k)))^{\text{T}}$.  
For an arbitrary vector norm $\Vert\cdot\Vert$, the induced norm of a matrix $W=(w_{1},...,w_{m})\in\mathbb{R}^{n\times m}$ is defined as $\Vert W\Vert=\sqrt{\sum_{i=1}^{m}\Vert w_{i}\Vert^{2}}$. Throughout, 
we slightly abuse the notation of vector norms and their induced matrix norms for simplicity. 

\vspace{-1ex}
\section{Preliminaries and Problem Formulation}

In this section, we first introduce some basics of graph theory. In what follows we formulate the problem of interest.
\vspace{-2ex}
\subsection{Basics of Graph Theory}
Consider a digraph $\mathcal{G}(\mathcal{V},\mathcal{E})$, 
where $\mathcal{V}=\{1,...,n\}$ denotes the set of nodes, $\mathcal{E}\subseteq \mathcal{V}\times \mathcal{V}$ represents the set of directed links, and $(i,j)\in\mathcal{E}$
implies that node $j$ can receive information from node $i$.
We denote $\mathcal{N}_{i}^{\text{in}}=\{j:(j,i)\in\mathcal{E}\}\cup \{i\}$ and
$\mathcal{N}_{i}^{\text{out}}=\{j:(i,j)\in\mathcal{E}\}\cup \{i\}$ as the in-neighbor set 
and out-neighbor set of node $i$, respectively.
A digraph is strongly connected if there exists a directed path between any pair of distinct nodes, 
which is commonly used in the literature \cite{NedicA2018Network}.
\begin{assumption} \label{Digraph}
The digraph $\mathcal{G}(\mathcal{V},\mathcal{E})$ is strongly connected.
\end{assumption}
\vspace{-2ex}
\subsection{Problem Formulation}

Consider the digraph $\mathcal{G}(\mathcal{V},\mathcal{E})$ where each node $i\in\mathcal{V}$ privately processes a convex function $f_{i}:\mathbb{R}^{m}\to\mathbb{R}$. 
All nodes collaboratively solve the following optimization problem\footnote{For clarity of presentation, we only consider the scalar variable case, i.e., $m=1$, as our algorithm and its analysis can be easily extended to the vector case by using the Kronecker operator.}:
\begin{eqnarray} \label{Problem}
\mathop{\text{minimize}}_{x\in\mathbb{R}^{m}} ~ f(x)=\sum_{i=1}^{n}f_{i}(x).
\end{eqnarray} 
In such a problem, each node $i$ maintains a local estimate $x_{i}(k)\in\mathbb{R}^{m}$ of the decision vector $x$ at each step $k$ and can only share its own information with a subset of nodes via the communication network. 
We make the following assumptions on the local functions:

\begin{assumption}\label{StrongConvex}
Each $f_{i}$ is $\mu$-strongly convex, i.e., there exists a $\mu>0$, such that 
$f_{i}(y)\geq f_{i}(x)+\nabla f_{i}(x)^{\text{T}}(y-x)+\frac{\mu}{2}\Vert x-y\Vert_{2}^{2}$, $\forall x,y\in\mathbb{R}^{m}$.
\end{assumption}


\begin{assumption}\label{LipschitzF}
Each $f_{i}$ is $L$-smooth, i.e., 
$\Vert \nabla f_{i}(x)-\nabla f_{i}(y)\Vert_{2}\leq L\Vert x-y\Vert_{2}$ for some $L>0$, $\forall x,y\in\mathbb{R}^{m}$.
\end{assumption}


Assumptions \ref{StrongConvex}-\ref{LipschitzF} are standard for 
the linear convergence in literature, see e.g. \cite{QuG2018Harnssing,XuJ2018Convergence}.
Under Assumption \ref{StrongConvex}, the problem (\ref{Problem}) has a unique optimal solution $x^{\star}\in\mathbb{R}^{m}$.


The main objective of this paper is to design a distributed algorithm where nodes are
only allowed to communicate quantized variables over $\mathcal{G}(\mathcal{V},\mathcal{E})$, with linear convergence to the exact optimal solution $x^{\star}$ of problem (\ref{Problem}).

\vspace{-1ex}
\section{Algorithm Development}\label{sec-3}

\begin{figure*}[htb]
\centering
\includegraphics[width=13cm]{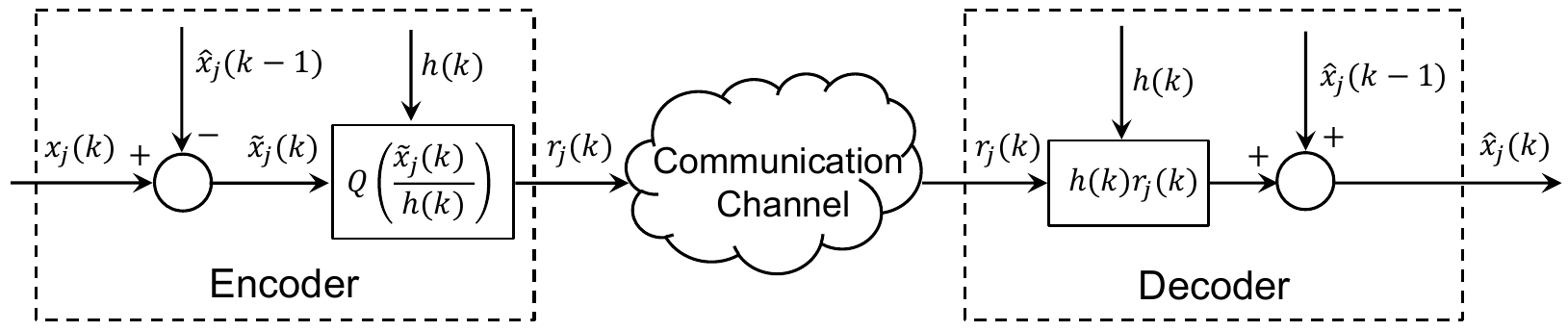}
\caption{The encoding-decoding scheme.}
\label{figStrcture}
\end{figure*}

In this section, we first explain why  a quantizer cannot be directly incorporated into 
the push-pull/$\mathcal{AB}$ Algorithm\cite{PuS2018PushJou,XinR2018A}. Then,  
 we propose the Q-DGT and show that it is robust to quantization errors. 
Finally, we introduce the 
quantization rule. 



%
%

\vspace{-2ex}
\subsection{Push-pull/$\mathcal{AB}$ Algorithm with Naive Quantization Does Not Work}

In the push-pull/$\mathcal{AB}$ algorithm\cite{PuS2018PushJou,XinR2018A}, each node $i\in\mathcal{V}$ maintains 
two vectors $x_{i}(k)$ and $y_{i}(k)$ per step $k\in\mathbb{N}$,  and performs the following updates: 
\begin{subequations} \label{MainAlg0DGT}
\begin{align}
x_{i}(k+1)
=\hspace*{0.5em}&
\sum_{j=1}^{n}a_{ij}x_{j}(k)-\eta y_{i}(k), \label{MainAlg0DGT-a}\\
y_{i}(k+1)
=\hspace*{0.5em}&
\sum_{j=1}^{n}b_{ij} y_{j}(k)+\nabla f_{i}(x_{i}(k+1)-\nabla f_{i}(x_{i}(k)), \label{MainAlg0DGT_2}
\end{align}
\end{subequations}
%
%
where 
$y_{i}(0)=\nabla f_{i}(x_{i}(0))$, 
$\mathcal{A}=[a_{ij}]_{n\times n}$ and $\mathcal{B}=[b_{ij}]_{n\times n}$ are weight matrices induced by $\mathcal{G}(\mathcal{V},\mathcal{E})$ 
satisfying: (1) $a_{ij}>0$ for $j\in\mathcal{N}_{i}^{\text{in}}$, otherwise $a_{ij}=0$, and $\sum_{j\in\mathcal{N}_{i}^{\text{in}}}a_{ij}=1$; (2) $b_{ij}>0$ for $i\in\mathcal{N}_{j}^{\text{out}}$, otherwise $b_{ij}=0$, and $\sum_{i\in\mathcal{N}_{j}^{\text{out}}}b_{ij}=1$.
Naive quantization means that $x_{j}(k)$ and $y_{j}(k)$ in (\ref{MainAlg0DGT}) are replaced by their quantized versions $\hat{x}_{j}(k)$ and $\hat{y}_{j}(k)$, respectively. That is, 
if $x_{j}(k)$ and $y_{j}(k)$ in (\ref{MainAlg0DGT}) are directly quantized as $\hat{x}_{j}(k)$ and $\hat{y}_{j}(k)$, respectively. 
Then, 
\begin{subequations} \label{MainAlg0DGTError}
\begin{align}
x_{i}(k+1)
=\hspace*{0.5em}&
\sum_{j=1}^{n}a_{ij}\hat{x}_{j}(k)-\eta y_{i}(k), \label{MainAlg0DGT-aError}\\
y_{i}(k+1)
=\hspace*{0.5em}&
\sum_{j=1}^{n}b_{ij} \hat{y}_{j}(k) 
+\nabla f_{i}(x_{i}(k+1)-\nabla f_{i}(x_{i}(k)), \label{MainAlg0DGT_2Error}
\end{align}
\end{subequations}
%
%
where $\sigma_{x_{j}}(k)\triangleq \hat{x}_{j}(k)-x_{j}(k)$ and $\sigma_{y_{j}}(k)\triangleq \hat{y}_{j}(k)-y_{j}(k)$ are the quantization errors.  
Note that if $\mathcal{A}=\mathcal{B}$, then (\ref{MainAlg0DGTError}) is exactly the quantized algorithm in \cite{Kajiyama2021Linear}. However, taking summation over $i\in\mathcal{V}$, (\ref{MainAlg0DGT_2Error}) implies that 
\begin{eqnarray}\label{MainAlg0DGTErrorAccu}
\bm{1}_{n}^{\text{T}}y(k+1)=\bm{1}_{n}^{\text{T}}\nabla F(x(k+1))+\sum_{l=0}^{k}\bm{1}_{n}^{\text{T}}\sigma_{y}(l).  
\end{eqnarray}
Thus the quantization errors are accumulated,  
i.e., $\sum_{l=0}^{k}\bm{1}_{n}^{\text{T}}\sigma_{y}(l)$. 
No matter whether $\sigma_{y}(k)$ in (\ref{MainAlg0DGTErrorAccu}) converges or not, 
 $\bm{1}_{n}^{\text{T}}y(k)$ cannot exactly track the global gradient $\bm{1}_{n}^{\text{T}}\nabla F(x(k))$. 

This observation was first pointed out in \cite{PuS2020AR} and then the author proposed a robust push-pull algorithm. 
However,  the work \cite{PuS2020AR} does not involve the design of quantizer and simply assumes that 
$\mathbb{E}[\sigma_{x}(k)]=\mathbb{E}[\sigma_{y}(k)]=\bm{0}_{n}$, and $\mathbb{E}[\Vert\sigma_{x}(k)\Vert^{2}]\leq\sigma_{x}$, $\mathbb{E}[\Vert\sigma_{y}(k)\Vert^{2}]\leq\sigma_{y}$ for some $\sigma_{x}, \sigma_{y}>0$. 
This condition is clearly not satisfied for the deterministic quantizers. 
In addition, the algorithm in \cite{PuS2020AR} can only converge to a neighborhood of the optimal solution in expectation. 
All above motivates us to propose the Q-DGT. 
%

\vspace{-1ex}
\subsection{The Q-DGT Algorithm}

In this work,
we design a dynamic encoding-decoding scheme for quantized communication (see Fig. 1). 
At step $k$, 
each node $j\in\mathcal{V}$ encodes $x_{j}(k)$ into 
$r_{j}(k)$  by using:
\begin{eqnarray}\label{rrrrjk}
r_{j}(k)
&\hspace*{-0.5em}=\hspace*{-0.5em}&
Q_{K_{x}}\left(\frac{1}{h(k)}\left(x_{j}(k)-\hat{x}_{j}(k-1)\right)\right), \label{encoder111}
\end{eqnarray}
where $\hat{x}_{j}(-1)=\bm{0}_{m}$ and 
$\hat{x}_{j}(k-1)$ is an estimation of $x_{j}(k-1)$, 
 $h(k)$ is a decaying scaling function. 
 Note that we quantize the scaled 
 ``innovation'', i.e., $\frac{1}{h(k)}(x_{j}(k)-\hat{x}_{j}(k-1))$. 
 The reason is that the amplitude of the prediction error is usually smaller than that of the state itself such that the scaled 
 	``innovation'' can be quantized by fewer bits. However, 
 	it brings challenge for the finite-level quantizer to avoid saturation. 
 We will show later that 
  the value of $x_{j}(k)-\hat{x}_{j}(k-1)$ decays to zero at the speed of the same order of $h(k)$, and 
rigorously prove that the scaled ``innovation'' can always be upper bounded by a finite constant.  
 Then, 
node $j$ broadcasts $r_{j}(k)$ to its out-neighbors. 
Upon $r_{j}(k)$ is received by 
the out-neighbor node $i\in\mathcal{N}_{j}^{\text{out}}$, 
it  
decodes $r_{j}(k)$  
as follows:
\begin{eqnarray}\label{hatttttxjk}
\hat{x}_{j}(k)&\hspace*{-0.5em}=\hspace*{-0.5em}&h(k)r_{j}(k)+\hat{x}_{j}(k-1).  \label{Decoder1} 
\end{eqnarray} 
Here $h(k)$ plays a critical role in estimating the states of node $j$. We highlight that all the out-neighbors of node $j$ receive the same information, so we do not distinguish the specific subscript.  
The above encoding-decoding scheme is performed for $y_{j}(k)$ in the same way, i.e., encode $y_{j}(k)$ into $s_{j}(k)$ and decode $s_{j}(k)$ to $\hat{y}_{j}(k)$.  
However,  
the deterministic quantization errors 
makes it infeasible to apply the robust push-pull algorithm \cite{PuS2020AR} directly in our setting (see Section III-A).  
To resolve it, we design the updates of node $i\in\mathcal{V}$ as follows: 
\begin{subequations} \label{MainAlg1}
\begin{align}
x_{i}(k+1)
=\hspace*{0.5em}&
x_{i}(k) +\alpha\sum_{j=1}^{n}a_{ij}\left(\hat{x}_{j}(k)-\hat{x}_{i}(k)\right)\notag\\
&-\eta \left(y_{i}(k)-y_{i}(k-1)\right), \label{MainAlg1-a}\\
y_{i}(k+1)
=\hspace*{0.5em}&
(1-\beta)y_{i}(k) +\beta\sum_{j=1}^{n}b_{ij}\hat{y}_{j}(k) \notag\\
&+\nabla f_{i}(x_{i}(k+1)), 
 \label{AlgResB_2}
\end{align}
\end{subequations}
where
$\alpha,\beta \in(0,1)$ are two positive constants,
 $\eta\geq 0$
is a constant step size
 that will be specified later.   
Although node $i\in\mathcal{V}$ can access to its true values $x_{i}(k)$ and $y_{i}(k)$ at step $k$, the estimate $\hat{x}_{i}(k)$ and $\hat{y}_{i}(k)$ are also used in our algorithm for error compensations, which is of the similar spirit as in the quantized average consensus in \cite{LiT2010Distributed}. 
 We summarize the Q-DGT in Algorithm \ref{algOnlineF}.

Now, we demonstrate why the Q-DGT  is robust to quantization errors. 
Let $\mathcal{A}_{\alpha}\triangleq (1-\alpha)I_{n}+\alpha\mathcal{A}$ and 
$ \mathcal{B}_{\beta} \triangleq (1-\beta)I_{n}+\beta\mathcal{B}$. 
Then, (\ref{MainAlg1}) can be rewritten as the following compact form: 
\begin{subequations} \label{MainAlgCom}
\begin{align}
x(k+1)
=\hspace*{0.5em}&
\mathcal{A}_{\alpha}x(k)+\alpha(\mathcal{A}-I_{n})\sigma_{x}(k)  \notag\\
&-\eta(y(k)- y(k-1)),  \label{AlgRes_1ARew}  \\
y(k+1)
=\hspace*{0.5em}&
\mathcal{B}_{\beta}y(k)+\nabla F(x(k+1))+\epsilon_{y}(k), \label{AlgRes_2BRew}
\end{align}
\end{subequations} 
where $\epsilon_{y_{i}}(k)\triangleq\beta \sum_{j=1}^{n}b_{ij}\sigma_{y_{j}}(k)$. 

Let $z(k)\triangleq y(k)- y(k-1)$.
Then, 
\begin{subequations} \label{MainAlgComRew}
\begin{align}
x(k+1)
=\hspace*{0.5em}&
\mathcal{A}_{\alpha}x(k)+\alpha(\mathcal{A}-I_{n})\sigma_{x}(k)-\eta z(k),  \label{AlgRes_1A}  \\
z(k+1)
=\hspace*{0.5em}&
\mathcal{B}_{\beta}z(k)+\left(\nabla F(x(k+1))+\epsilon_{y}(k)\right)  \notag\\
&-\left(\nabla F(x(k))+\epsilon_{y}(k-1)\right). \label{AlgRes_2B}
\end{align}
\end{subequations}
Assumption \ref{Digraph} implies that $\mathcal{A}$ has a unique nonnegative left eigenvector $\pi_{\mathcal{A}}$ such that $\pi_{\mathcal{A}}^{\text{T}}\bm{1}_{n}=1$ and $\pi_{\mathcal{A}}^{\text{T}}\mathcal{A}=\pi_{\mathcal{A}}^{\text{T}}$, and $\mathcal{B}$ has a unique nonnegative right eigenvector $\pi_{\mathcal{B}}$ such that $\pi_{\mathcal{B}}^{\text{T}}\bm{1}_{n}=1$ and $\mathcal{B}\pi_{\mathcal{B}}=\pi_{\mathcal{B}}$ \cite{ZhangJ2020Distributed}.   
%
%
%
Define $\bar{x}(k)\triangleq \pi_{\mathcal{A}}^{\text{T}}x(k)$ and $\bar{z}(k)\triangleq\bm{1}_{n}^{\text{T}}z(k)$, 
we obtain 
\begin{subequations} 
	\begin{align}
		\bar{x}(k+1)
		=\hspace*{0.5em}&
		\bar{x}(k)-\eta \pi_{\mathcal{A}}^{\text{T}}z(k), \label{barxx} \\
		\bar{z}(k+1)
		=\hspace*{0.5em}&
	\bar{z}(k)+\bm{1}_{n}^{\text{T}}\left(\nabla F(x(k+1))+\epsilon_{y}(k)\right)  \notag\\
		&-\bm{1}_{n}^{\text{T}}\left(\nabla F(x(k))+\epsilon_{y}(k-1)\right).  \label{yyyya}
	\end{align}
\end{subequations}
%
%
%
%
%
%
%
Conducting mathematical induction for (\ref{yyyya}) yields that 
\begin{eqnarray} \label{Important111}
\bar{z}(k+1)
&\hspace*{-0.5em}=\hspace*{-0.5em}&
\bm{1}_{n}^{\text{T}}\left(\nabla F(x(k+1))+\epsilon_{y}(k)\right). 
\end{eqnarray}
Notably, 
 the accumulated error $\sum_{l=0}^{k}\bm{1}_{n}^{\text{T}}\sigma_{y}(l)$ in (\ref{MainAlg0DGTErrorAccu}) disappears in (\ref{Important111}).  
If $\bm{1}_{n}^{\text{T}}\epsilon_{y}(k)$ tends to zero, then $\bar{z}(k+1)$ tends to the exact global gradient $\bm{1}_{n}^{\text{T}}\nabla F(x(k+1))$. 
In contrast to \cite{PuS2020AR}, we do not make any assumption on the error $\epsilon_{y}(k)$.  
This requires to design the Q-DGT (\ref{MainAlg1}) carefully and 
handle the joint effects of quantization errors on $x(k)$ and $z(k)$. 
Specifically, our algorithm can converge linearly and even support 3-level quantization.
\begin{algorithm}[t]
\caption{The Q-DGT{\color{black}~~---from the view of node $i$}}
\label{algOnlineF}
\begin{algorithmic}[1]
    \State \textbf{Initialization:} randomly initialize $x_{i,0}$, and $y_{i,0}$ for each $i\in\mathcal{V}$.

    \For {$k=0,1,2,...$}

    \State\textbf{Encoder:} calculate $r_{i}(k)$ and $s_{i}(k)$.

   \State\textbf{Communication:}  broadcast $r_{i}(k)$ and $s_{i}(k)$ to its out-neighbors, 
 and  receive $r_{j}(k)$ and  $s_{j}(k)$ from its in-neighbors $j\in\mathcal{N}_{i}^{\text{in}}$.


\State \textbf{Decoder:} calculate $\hat{x}_{j}(k)$ and $\hat{y}_{j}(k)$.

\State \textbf{Updation:} update $x_{i}(k+1)$ and $y_{i}(k+1)$ via (\ref{MainAlg1}).
%
%
    \EndFor

    \State \textbf{Return}: $\{x_{i}(k)\}$.

\end{algorithmic}
\end{algorithm}

\vspace{-1ex}
\subsection{The Quantization Rule}
The uniform quantizer $Q_{K}(\cdot)$ for a vector $u=(u_{1},...,u_{m})^{\text{T}}$ is defined as 
$Q_{K}(u)=\left(q(u_{1}),...,q(u_{m})\right)^{\text{T}}$ with 
\begin{equation}\label{Quan2}
q(u_{i})=
\begin{cases}
0, &\mbox{$-1/2<u_{i}\leq1/2$}\\
k, &\mbox{$\frac{2 k-1}{2}< u_{i}\leq\frac{2 k+1}{2}$, $k=1,...,K$}\\
K, &\mbox{$u_{i}>\frac{2K+1}{2}$}\\
-q(-u_{i}), &\mbox{$u_{i}\leqslant -1/2$}
\end{cases}\notag
\end{equation}
for $i=1,...,m$. 
The quantizer $q(\cdot)$ maps a real number to a finite set 
$\mathcal{S}=\{0,\pm k; k=1,2,...,K\}$ with $K\in\mathbb{N}_{+}$. 
The quantization level of $q(\cdot)$ is $2K+1$.  
If $\|u\|_{\infty}\leq K+1/2$, 
 the quantizer is not saturated, and the quantization error is bounded, i.e., $\|u-Q_{K}(u)\|_{\infty}\leqslant1/2$.

\section{Convergence Analysis}\label{sec-IV}

In this section, 
we first establish 
lower bounds for the quantization levels  
to solve the saturation issue. 
Then, 
the linear convergence of Q-DGT under finite-level quantization  is rigorously proved. 
Finally, we show that Q-DGT converges linearly even 
with $3$-level quantization. 
%

\vspace{-1ex}
\subsection{Design of Finite Quantization Levels to Avoid Saturation}

Note that the joint effect of quantization on the evolutions of $x(k)$ and $z(k)$ brings challenges to design 
the finite quantization levels.  
To solve this issue, 
we first derive the upper bound of the feasible step size, 
and then obtain the lower bounds for the quantization levels. 

\begin{lemma}[\cite{ZhangJ2020Distributed,PuS2018PushJou}]\label{Normequ2}
Suppose Assumption \ref{Digraph} holds. There exists matrix norms $\Vert\cdot\Vert_{\mathcal{A}}$
and $\Vert\cdot\Vert_{\mathcal{B}}$ such that
$\sigma_{\mathcal{A}}\triangleq \Vert\mathcal{A}_{\alpha}-\bm{1}_{n}\pi_{\mathcal{A}}^{\text{T}}\Vert_{\mathcal{A}}<1$
and
$\sigma_{\mathcal{B}}\triangleq \Vert\mathcal{B}_{\beta}-\pi_{\mathcal{B}}\bm{1}_{n}^{\text{T}}\Vert_{\mathcal{B}}<1$.
Moreover, there exists positive scalars $\delta_{\mathcal{A}2}$,
$\delta_{\mathcal{B}2}$, $\delta_{\mathcal{AB}}$ and $\delta_{\mathcal{BA}}$ such that for any
$X\in\mathbb{R}^{n\times p}$, we have 
$\delta_{\mathcal{BA}}^{-1}\Vert X \Vert_{\mathcal{B}}\leq\Vert X \Vert_{\mathcal{A}}\leq \delta_{\mathcal{AB}}\Vert X \Vert_{\mathcal{B}}$, 
$\Vert X \Vert_{2}\leq\Vert X \Vert_{\mathcal{B}}\leq \delta_{\mathcal{B}2}\Vert X \Vert_{2}$, 
and 
$\delta_{\mathcal{A}2}^{-1}\Vert X \Vert_{\mathcal{A}}\leq \Vert X \Vert_{2}\leq\Vert X \Vert_{\mathcal{A}}$.
\end{lemma}
Define 
\begin{align}\label{Thetakkkkk}
	\Theta(k)\triangleq \left(\Vert \bar{x}(k)-x^{\star}\Vert_{2},\Vert x(k)-\bm{1}_{n}\bar{x}(k)\Vert_{\mathcal{A}},
	\Vert z(k)-\pi_{\mathcal{B}}\bar{z}(k)\Vert_{\mathcal{B}}\right)^{\text{T}}. 
\end{align}
%
To facilitate the subsequent analysis, 
we further define: 
$\kappa_{1}\triangleq \Vert I_{n}-\bm{1}_{n}\pi_{\mathcal{A}}^{\text{T}}\Vert_{\mathcal{A}}$,
$\kappa_{2}\triangleq \Vert \pi_{\mathcal{B}}\Vert_{\mathcal{A}}$,
$\kappa_{3}\triangleq \Vert I_{n}-\pi_{\mathcal{B}}\bm{1}_{n}^{\text{T}}\Vert_{\mathcal{B}}$,
$\kappa_{4}\triangleq \Vert \mathcal{A}_{\alpha}-I_{n}\Vert_{2}$. 
The following lemma provides a linear matrix inequality, which
will be instrumental in establishing 
the lower bound for quantization level. 
\begin{lemma}\label{proposition1}
Suppose Assumptions 1-3 hold. If 
 the step size $\eta\leq \frac{1}{(\mu+L)\pi_{\mathcal{A}}^{\text{T}}\pi_{\mathcal{B}}}$,   
then 
\begin{eqnarray}\label{System}
\Theta(k+1)\preceq G\Theta(k)+\varsigma(k),
\end{eqnarray}
where 
$\Theta(k)$ is defined in (\ref{Thetakkkkk}),  
 the notation $\preceq$ means element-wise less than or equal to,
$G\in\mathbb{R}^{3\times 3}$ and $\varsigma(k)\in\mathbb{R}^{3\times 1}$ are given by (\ref{GmatrixGG}) and (\ref{varsigma123}), respectively. 
\end{lemma}

\textit{Proof:}
See Appendix A. 
$\hfill \blacksquare$

In Lemma \ref{proposition1}, the presence of $\varsigma(k)$ is due to quantization errors.  
If $\varsigma(k)$ linearly converges to $\bm{0}_{3}$, 
then we can prove that $\Theta(k)$ linearly converges to $\bm{0}_{3}$ provided that the spectral radius $\rho(G)<1$. After that, the linear convergence of Q-DGT can be proved.
%
%
We first provide a sufficient condition in terms of the step size $\eta$ to guarantee $\rho(G)< 1$.

\begin{lemma}\label{StepCond}
Suppose Assumptions 1-3 hold.
If the step size $\eta$ satisfies 
%
	\begin{align}\label{Range}
		\eta \leq \min\Bigg\{\frac{1}{(\mu+L)\pi_{\mathcal{A}}^{\text{T}}\pi_{\mathcal{B}}},\frac{1-\sigma_{\mathcal{A}}}{2\sqrt{n}\kappa_1\kappa_2 L\delta_{\mathcal{A}2}},
		\frac{1-\sigma_{\mathcal{B}}}{2\delta_{\mathcal{B}2}\kappa_3 L}, ~~~\notag\\  \frac{2\Gamma_3}{\Gamma_2+\sqrt{\Gamma_2^2+4\Gamma_1\Gamma_3}}\Bigg\},   
	\end{align}
where $\Gamma_{i}$, $i=1,2,3$, are constants given by (\ref{GXi1}). Then,
$\rho(G)<1$.  
\end{lemma}

\textit{Proof:}
See Appendix B. 
$\hfill \blacksquare$

\begin{remark}
Note that the network information is required to calculate the upper bound of step size. 
If $\mathcal{A}$ and $\mathcal{B}$ are known, then the parameters $\kappa_{1}$-$\kappa_{4}$, 
the step size $\eta$, and the spectral radius $\rho(G)$ can be obtained. 
Typically, this requirement is necessary even for unquantized push-pull/$\mathcal{AB}$ algorithms\cite{PuS2018PushJou,XinR2018A,PuS2020AR}.
\end{remark}

%
%

It is known that if $\rho(G)<1$,
then there exist a matrix norm $\|\cdot\|_{\mathcal{C}}$ and a constant $\tau$ such that 
$\|G^k\|_{2} \leq  \tau\hat{\rho}^k$
for an arbitrarily small  constant $\varpi>0$\cite{HornR2012Matrix},  where  
\begin{eqnarray}
\hat{\rho}\triangleq\rho(G)+\varpi<1.
\end{eqnarray}

Now, we are in a position to provide conditions on the quantization levels, under which 
	the saturation issue can be solved. 



\begin{theorem}\label{Unsat1}
Suppose Assumptions 1-3 hold. Let $h(k)=C\xi^{k}$,
where $C$ is a known positive constant, and $\xi\in(\hat{\rho},1)$.
The step size $\eta$ is chosen according to (\ref{Range}).
Then the quantizers will never saturate provided that 
$K_{x}$ and $K_{y}$ satisfy the following conditions: 
\begin{eqnarray}\label{QuanLevelConAlterQQQ}
\hspace*{-1em}	K_x &\hspace*{-0.5em}\geq\hspace*{-0.5em}& 
	\max\Bigg\{\frac{v_{1}}{C} - \frac{1}{2}, ~
	\frac{\sqrt{3}\varphi_{1}\Vert\Theta(0)\Vert_{2}}{C \xi} + \frac{2\alpha n+1}{2\xi}  - \frac{1}{2},~  \notag\\
\hspace*{-1em}&&~~~~~~~~~	\frac{\sqrt{3}\varphi_{1} \tau \Vert\Theta(0)\Vert_{2}}{C \xi} \bar{\Upsilon} + \frac{2\alpha n+1}{2\xi} + \frac{ n\eta\beta }{2  \xi^2} - \frac{1}{2}\Bigg\}, \notag\\
\hspace*{-1em}	K_y &\hspace*{-0.5em}\geq\hspace*{-0.5em}&
	 \max\Bigg\{\frac{v_{2}}{C} - \frac{1}{2}, ~
	\frac{\sqrt{3}\varphi_{2}\tau \Vert\Theta(0)\Vert_{2}}{C} \bar{\Upsilon} + \frac{n \beta +1 }{2\xi} - \frac{1}{2}\Bigg\}, 
\end{eqnarray}
where 
$v_{1}\triangleq\max_{i\in\mathcal{V}}\Vert x_{i}(0)\Vert_{\infty}$,  
$v_{2}\triangleq\max_{i\in\mathcal{V}}\Vert y_{i}(0)\Vert_{\infty}$, 
$\varphi_{1}$ and $\varphi_{2}$ are given in (\ref{varpi123}), 
$\bar{\Upsilon}\triangleq 1+\frac{\tilde{\varsigma}\hat{\rho}}{\xi(\xi-\hat{\rho})\Vert\Theta(0)\Vert_{2}}+\frac{\tilde{\varsigma}}{\xi\tau \Vert\Theta(0)\Vert_{2}}$ with 
the constant $\tilde{\varsigma}$ given by (\ref{tildevarsi}). 
\end{theorem}

\textit{Proof:} See Appendix D. 
$\hfill \blacksquare$ 

\begin{remark}\label{ruleExp} 
Theorem \ref{Unsat1} provides a sufficient condition to guarantee that the quantizers will never saturate.  
Note that all the terms on the right sides of (\ref{QuanLevelConAlterQQQ}) are finite constants, which 
	implies that the quantizers will never saturate as long as $K_{x}$ and $K_{y}$ are positive integers larger than the lower bounds in (\ref{QuanLevelConAlterQQQ}). 
In addition, (\ref{QuanLevelConAlterQQQ}) depends on the initial states of nodes, which is common in literature \cite{LiH2020Computation,YiP2014Quantized}. When executing the proposed algorithm in practice, 
we can choose $\alpha$ and $\beta$ from $(0,1)$ arbitrarily, let $\xi$ be in close proximity to $1$, and select a large enough constant and a small enough constant as the quantization level and the step size, respectively. 
\end{remark}

\vspace{-1ex}
\subsection{Linear Convergence under Finite Quantization Levels}

Building upon the conditions on the quantization levels in Theorem \ref{Unsat1}, 
the following theorem shows that the Q-DGT can linearly converge to the optimal solution at the rate of $\mathcal{O}(\xi^{k})$ with $\xi\in (\hat{\rho}, 1)$. 
\begin{theorem}
Suppose the conditions in Theorem \ref{Unsat1} are satisfied. Let $\{x_{i}(k)\}$, $i\in\mathcal{V}$, be the sequence generated by Algorithm \ref{algOnlineF}. If the quantization levels satisfy (\ref{QuanLevelConAlterQQQ}), then  
Q-DGT can linearly converge to $x^{\star}$
at the rate of $\mathcal{O}(\xi^{k})$, i.e., $\Vert x_{i}(k)-x^{\star}\Vert_{2}=\mathcal{O}(\xi^{k})$ for all $i\in\mathcal{V}$.
%
\end{theorem}

\textit{Proof:}
Recalling (\ref{System}) and (\ref{tildeHl}), we can straightly obtain 
\begin{eqnarray}\label{ConverIneq}
\Vert \Theta(k)\Vert_{2}
&\hspace*{-0.5em}\leq\hspace*{-0.5em}&
 \Vert G^{k}\Vert_{2} \Vert \Theta(0)\Vert_{2}+\tilde{\varsigma}\sum_{l=0}^{k-1}\Vert G^{k-1-l}\Vert_{2}\xi^{l}\notag\\
 &\hspace*{-0.5em}\leq\hspace*{-0.5em}&
 \tau \hat{\rho}^{k} \Vert \Theta (0)\Vert_{2}+\tilde{\varsigma}\tau\sum_{l=0}^{k-1} \hat{\rho}^{k-1-l}\xi^{l} \notag
\end{eqnarray}
Note that 
$\sum_{l=0}^{k-1} \hat{\rho}^{k-1-l}\xi^{l}=\xi^{k-1}\sum_{l=0}^{k-1}
\left(\frac{\hat{\rho}}{\xi}\right)^{k-1-l}
\leq\frac{\xi^{k}}{\xi-\hat{\rho}}$. 
%
Hence, 
\begin{eqnarray}\label{ConverLinearR}
\Vert \Theta(k)\Vert_{2}
\leq
 \tau d_{0} \hat{\rho}^{k}
 +\frac{\tilde{\varsigma}\tau}{\xi-\hat{\rho}}\xi^{k}, \notag
\end{eqnarray}
which implies that $\Vert \Theta(k)\Vert_{2}$ converges to $0$ at the rate of $\mathcal{O}(\xi^{k})$.
Therefore, $\Vert \bar{x}(k)-x^{\star}\Vert_{2}$,
$\Vert x(k)-\bm{1}_{n}\bar{x}(k)\Vert_{\mathcal{A}}$
and $ \Vert z(k)-\pi_{\mathcal{B}}\bar{z}(k)\Vert_{\mathcal{B}}$
all linearly converge to $0$ 
at the same rate. 
Note that
\begin{eqnarray}\label{}
\Vert x_{i}(k)-x^{\star}\Vert_{2} 
&\hspace*{-0.5em}\leq\hspace*{-0.5em}&
\Vert x_{i}(k)-\bar{x}(k)\Vert_{2}+\Vert \bar{x}(k)-x^{\star}\Vert_{2} \notag\\
&\hspace*{-0.5em}\leq\hspace*{-0.5em}&
\Vert x(k)-\bm{1}_{n}\bar{x}(k)\Vert_{2}+\Vert \bar{x}(k)-x^{\star}\Vert_{2}, \notag
\end{eqnarray}
which implies that $\Vert x_{i}(k)-x^{\star}\Vert_{2}=\mathcal{O}(\xi^{k})$ for all $i\in\mathcal{V}$
by recalling the fact that $\Vert x(k)-\bm{1}_{n}\bar{x}(k)\Vert_{2}\leq \Vert x(k)-\bm{1}_{n}\bar{x}(k)\Vert_{\mathcal{A}}$.
\hfill $\blacksquare$

%


\vspace{-1ex}
\subsection{3-Level Quantization is Enough for Linear Convergence}
As shown in Theorem \ref{Unsat1}, the lower bounds in (\ref{QuanLevelConAlterQQQ}) are finite.  
This inspires us to consider whether there exists a minimum number of quantization level that can preserve the linear  convergence? 
The following theorem gives a positive answer and reveals that we can set $K_{x}=K_{y}=1$ by appropriately tuning the associated parameters. 
In such an extreme scenario, each node $i\in\mathcal{V}$ can solve problem (\ref{Problem}) with 
	$3$-level quantization. 

\begin{theorem}\label{MinimalQLThe}
Suppose the conditions in Theorem \ref{Unsat1} are satisfied.
If $\alpha$ and $ \beta$ are 
sufficiently small,
then there exists $C>0$ and $\xi\in (\hat{\rho},1)$ such that 
$K_{x}=K_{y}=1$ is sufficient to guarantee the linear convergence of Q-DGT. 
\end{theorem}

\textit{Proof:} 
To prove the result, our strategy is 
minimizing the lower bounds obtained in (\ref{QuanLevelConAlterQQQ}).   
Particularly, if we can choose the associated parameters appropriately such that all the lower bounds in (\ref{QuanLevelConAlterQQQ}) can be upper bounded by 1, then it can be concluded that the quantizers will never saturate even when $K_{x}=K_{y}=1$. In this case, Theorem 3 can be proved by recalling Theorem 2. 

Now, we consider the last term of each inequality in (\ref{QuanLevelConAlterQQQ}). 
Recalling 
$\bar{\Upsilon}\triangleq  1+\frac{\tilde{\varsigma}\hat{\rho}}{\xi(\xi-\hat{\rho})\Vert\Theta(0)\Vert_{2}}+\frac{\tilde{\varsigma}}{\xi\tau \Vert\Theta(0)\Vert_{2}}$ in Theorem \ref{Unsat1} 
and 
the expression of $\tilde{\varsigma}$ in (\ref{tildevarsi}).  
If $\alpha$ and $\beta$ both tend to 0, then $\tilde{\varsigma}$ tends to 0. 
Since $\Vert \Theta(0) \Vert_{2}$, $\hat{\rho}$ and $\tau$ are all some positive constants, and $\xi$ is a constant chosen in the interval $(\hat{\rho},1)$,  
we can obtain that $\bar{\Upsilon}$ tends to $1$. 
Therefore, the last term of each inequality in (\ref{QuanLevelConAlterQQQ}) 
can be upper bounded by $\frac{\sqrt{3}\varphi_{1} \tau \Vert\Theta(0)\Vert_{2}}{C\xi}+\frac{1}{2\xi}$ and $\frac{\sqrt{3}\varphi_{2} \tau \Vert\Theta(0)\Vert_{2}}{C}+\frac{1}{2\xi}$, respectively. 
Note that $\varphi_{1}$ and $\varphi_{2}$ are constants given in (\ref{varpi123}). 
If we choose the constant $\xi\in(\max\{0.5,\hat{\rho}\},1)$, and set 
$C>\max\left\{\frac{2\sqrt{3}\varphi_{1}\tau \Vert \Theta(0) \Vert_{2}}{2\xi-1}, \frac{2\sqrt{3}\varphi_{2} \xi\tau \Vert \Theta(0) \Vert_{2}}{2\xi-1}\right\}$, 
then the last term of each inequality in (\ref{QuanLevelConAlterQQQ}) both can be upper bounded by $1$.  
In addition,   
the other terms on the right side of (\ref{QuanLevelConAlterQQQ}) can be upper bounded by $1$ directly if we set 
$C>\max\{\frac{2}{3}v_{1}, \frac{2}{3}v_{2}, \frac{2\sqrt{3}\varphi_{1}\Vert \Theta(0) \Vert_{2}}{2\xi-1}\}$. 
In summary, there exists constants 
$\xi\in(\max\{0.5,\hat{\rho}\},1)$, and 
\begin{eqnarray} 
	C>\max\Bigg\{\frac{2}{3}v_{1}, \frac{2}{3}v_{2}, \frac{2\sqrt{3}\varphi_{1}\Vert \Theta(0) \Vert_{2}}{2\xi-1}, \frac{2\sqrt{3}\varphi_{1}\tau \Vert \Theta(0) \Vert_{2}}{2\xi-1}, \notag\\
	\frac{2\sqrt{3}\varphi_{2} \xi\tau \Vert \Theta(0) \Vert_{2}}{2\xi-1}\Bigg\}
\end{eqnarray} 
such that $K_{x}=K_{y}=1$  
is sufficient to guarantee the linear convergence of Q-DGT. 
\hfill $\blacksquare$

\begin{remark}
Quantized distributed algorithms with $3$-level quantization have also been studied in \cite{LiT2010Distributed,LeiJ2020Distributed,Zhang2019Distributed,YiP2014Quantized} to improve the communication efficiency. 
In contrast to the distributed optimization algorithms in \cite{LeiJ2020Distributed,Zhang2019Distributed,YiP2014Quantized}, the proposed Q-DGT can achieve linear convergence. 
Though the quantizer has only $3$ quantization levels, each node $i$ can still estimate the values of its in-neighbors $j\in\mathcal{N}_{i}^{\text{in}}$ iteratively via the decoding scheme (\ref{hatttttxjk}). 
This can be observed from the facts that $\hat{x}_{j}(k)-x_{j}(k)=h(k)e_{x_{j}}(k)$ and $\hat{y}_{j}(k)-y_{j}(k)=h(k)e_{y_{j}}(k)$. 
If the quantizers never saturate, then the diminishing $h(k)$ guarantees that $\hat{x}_{j}(k)$ and $\hat{y}_{j}(k)$ tend to $x_{j}(k)$ and $y_{j}(k)$, respectively, as $k$ tends to infinity. That is why our algorithm can converge to the true solution even with $3$-level quantization. 
%
%
\end{remark}




\section{Numerical Examples}\label{sec-V}

\begin{figure}
  \centering
  \includegraphics[width=4.2cm]{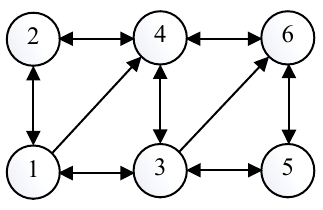}
  \caption{The directed graph.}\label{fig1}
\end{figure}

In this section, we apply our algorithm to the sensor fusion problem in
directed networks, which has been widely adopted in the literature \cite{PuS2018PushJou,XuJ2018Convergence}.
In this problem, all sensors collectively solve the following optimization problem over the digraph decipted in Fig. \ref{fig1}:
\begin{eqnarray} \label{SimulationPro}
\mathop{\text{minimize}}_{x\in\mathbb{R}^{m}}
~ f(x)=\sum_{i=1}^{n}\left(\Vert \mathcal{M}_{i}x-\zeta_{i}\Vert^{2}+\frac{\lambda}{2n}\Vert x\Vert^{2}\right), \notag
\end{eqnarray}
where 
$\mathcal{M}_{i}\in\mathbb{R}^{s\times m}$ and $\zeta_{i}\in\mathbb{R}^{s}$ denote the measurement matrix and the noise observation of sensor $i$, respectively,  $\lambda>0$ is the regularization parameter.

\begin{figure}[htbp]
\centering
\begin{minipage}{0.23\textwidth}
  \includegraphics[width=1\textwidth]{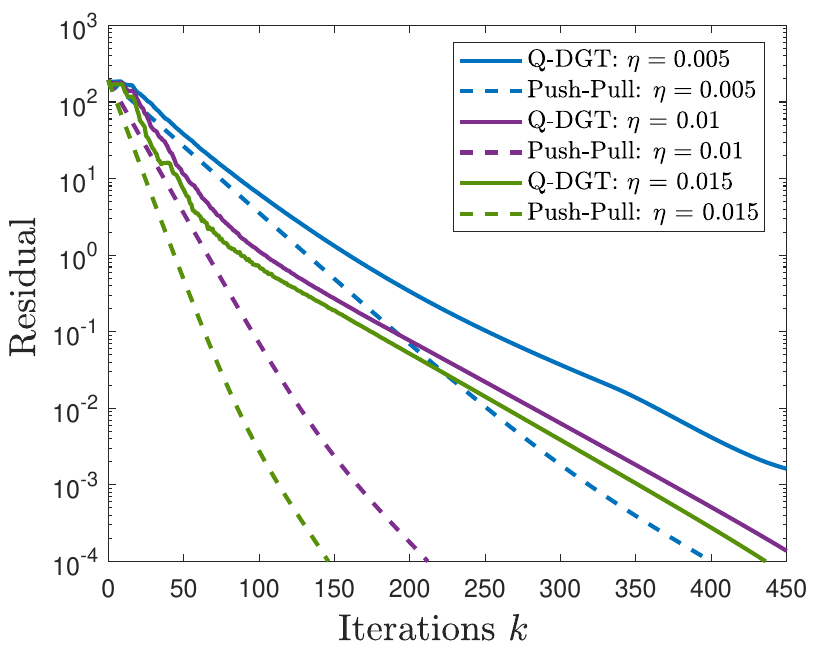}
  \caption*{(a)}
  \label{fig3-a}
\end{minipage}
\begin{minipage}{0.23\textwidth}
  \includegraphics[width=1\textwidth]{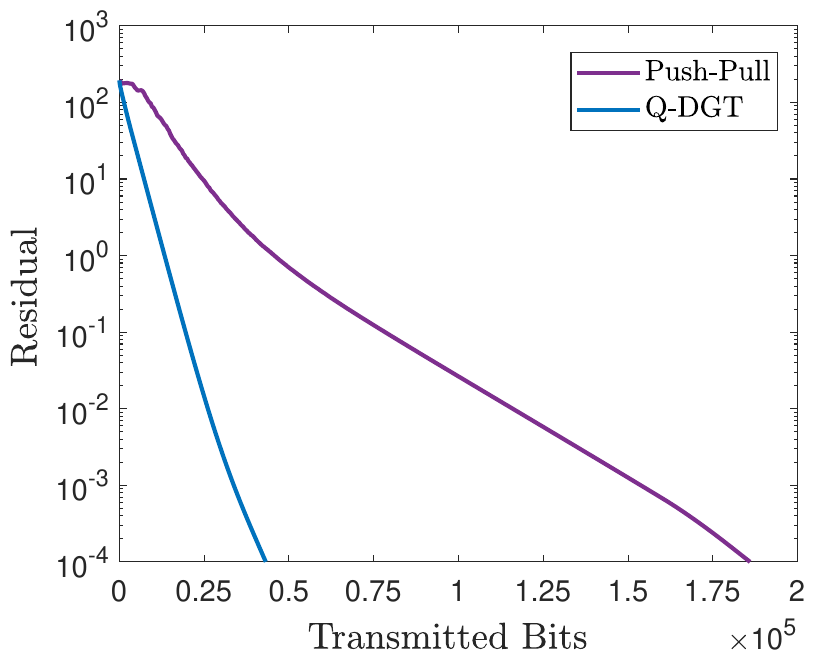}
  \caption*{(b)}
  \label{fig3-b}
\end{minipage}
\caption{(a) Comparison with push-pull in \cite{PuS2018PushJou} under different step sizes; 
(b) The total communication cost of Q-DGT and push-pull.}
\label{fig3}
\end{figure}

\begin{figure}[htb]
\centering
\begin{minipage}{0.23\textwidth}
  \includegraphics[width=1\textwidth]{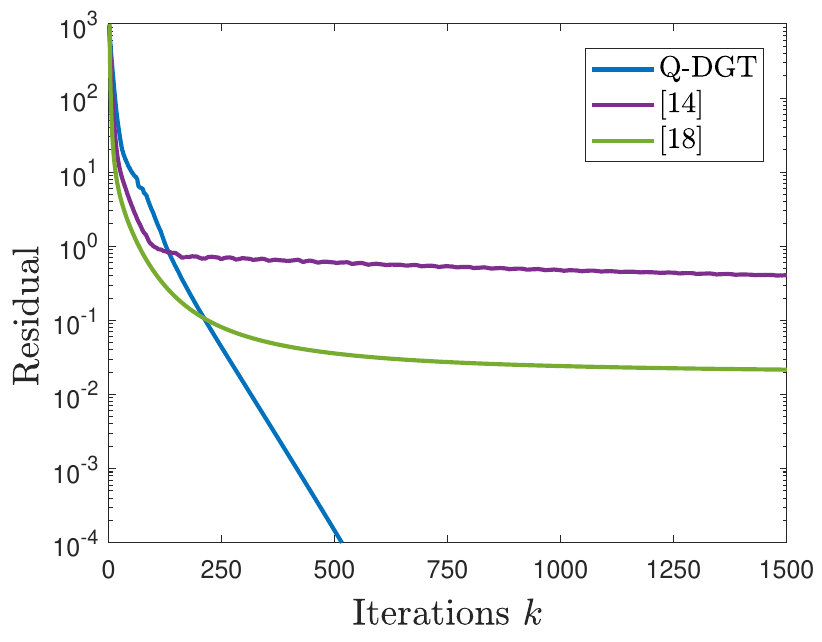}
  \caption*{(a)}
  \label{fig3-b}
\end{minipage}
\begin{minipage}{0.23\textwidth}
  \includegraphics[width=1\textwidth]{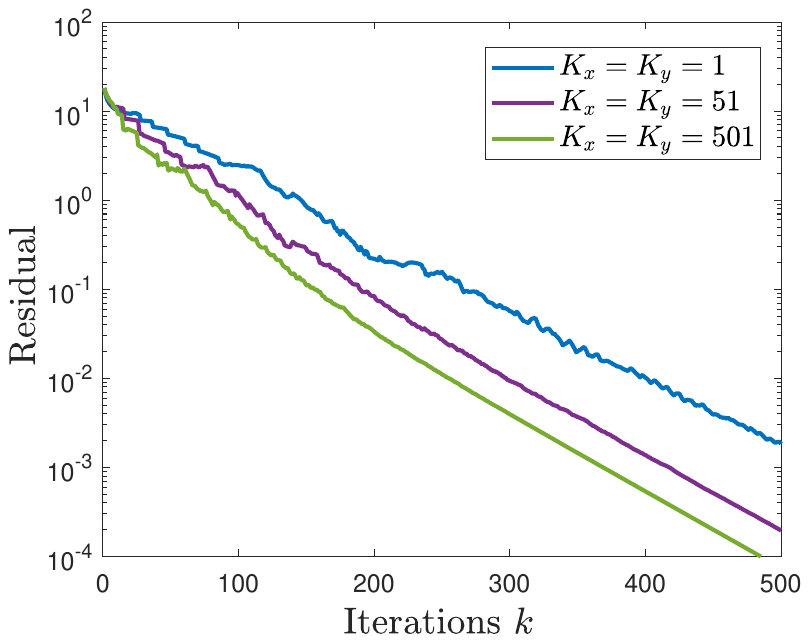}
  \caption*{(b)}
  \label{fig3-c}
\end{minipage}
\caption{(a) Comparisons with the quantized algorithms in \cite{YiP2014Quantized} and \cite{Taheri2020Quantized}; (b) Performances of Q-DGT under different  quantization levels.}
\label{fig4}
\end{figure}

In our simulations,  $\mathcal{M}_{i}\in\mathbb{R}^{2\times 2}$ and $\zeta_{i}\in\mathbb{R}^{2}$ are generated randomly for each $i\in\mathcal{V}$. We set $\lambda=0.05$. 
$\mathcal{A}$ and $\mathcal{B}$ are designed according to the rules in Remark 2 of \cite{PuS2018PushJou}.  
We first compare the convergence performance of Q-DGT with push-pull algorithm \cite{PuS2018PushJou} under different stepsizes.  The simulation results are depicted in Fig. \ref{fig3}(a). We can find that Q-DGT converges slower than the push-pull algorithm, which is reasonable as the performance inevitably affected by the loss of information. Despite this, the Q-DGT still maintains linear convergence, which is consistent with our theoretical results.
We further compare the total cost of communicated bits between the two algorithms with $\eta=0.008$. As shown in Fig. \ref{fig3}(b), the proposed Q-DGT requires less communicated bits for achieving the equal accuracy.  
Then, we make comparisons with the subgradient-based quantized distributed algorithms in \cite{YiP2014Quantized} and \cite{Taheri2020Quantized}. For fair comparison, we neglect the directionality in Fig. \ref{fig1} and adopt 
(\ref{QuanLevelConAlterQQQ}) for Q-DGT.
The results are depicted in Fig. \ref{fig4}(a).  It can be seen that the convergence rate of Q-DGT outperforms that of the quantized algorithms in \cite{YiP2014Quantized} and \cite{Taheri2020Quantized}.
Finally, we verify the effectiveness of Q-DGT under different fixed numbers of quantization levels.
The related parameters are chosen heuristically to meet the requirements in Theorem \ref{MinimalQLThe}.
As we can see in Fig. \ref{fig4}(b), the Q-DGT can still achieve linear convergence, 
even when the exchanged variables are respectively 
quantized with 3 quantization levels. 
In addition, a larger quantization level leads to faster convergence. This result is also reasonable since a larger quantization level implies a smaller quantization error.


%

\vspace{-1ex}

\section{Conclusion}\label{sec-V}

In this paper, we have studied the distributed optimization problem over directed networks with 
quantized communications. 
To cope with this problem, 
a novel quantized distributed algorithm Q-DGT has been proposed. 
The lower bounds for the number of quantization levels have been explicitly derived.  
We have rigorously shown that 
Q-DGT is robust to quantization errors, 
and achieves linear convergence 
even when the exchanged variables are respectively 
quantized with 3 quantization levels.  
Future works can focus on extending the proposed algorithm to time-varying directed networks.
It is also of interest to relax the conditions that preserves the convergence performance.

\vspace{-1ex}

\section*{Appendix }

\subsection{\textit{Proof of Lemma \ref{proposition1}.}}
For the clarity of presentation, we define 
$g(k)\triangleq \bm{1}_{n}^{\text{T}}\nabla F(x(k))$ and $\bar{g}(k)\triangleq\bm{1}_{n}^{\text{T}}\nabla F(\bm{1}_{n}\bar{x}(k))$. To prove this lemma, we first provide the following intermediate result. 
%
\begin{lemma} \label{Found1}
Suppose Assumptions \ref{StrongConvex}-\ref{LipschitzF} hold. We have
$\Vert g(k)-  \bar{g}(k)\Vert_{2}\leq \sqrt{n}L\Vert x(k)-\bm{1}_{n}\bar{x}(k)\Vert_{2}$,
$\Vert \bar{z}(k)-g(k)\Vert_{2}\leq \Vert  \bm{1}_{n}^{\text{T}}\epsilon_{y}(k-1)  \Vert_{2}$,
and $\Vert \bar{g}(k)\Vert_{2}\leq nL\Vert \bar{x}(k)-x^{\star}\Vert_{2}$.
If $\eta\leq \frac{1}{(\mu+L)\pi_{\mathcal{A}}^{\text{T}}\pi_{\mathcal{B}}}$,
then   
\begin{eqnarray}
\Vert \bar{x}(k)-\eta \pi_{\mathcal{A}}^{\text{T}}\pi_{\mathcal{B}}\bar{g}(k)-x^{\star}\Vert_{2}
&\hspace*{-0.5em}\leq\hspace*{-0.5em}&
(1-\eta \pi_{\mathcal{A}}^{\text{T}}\pi_{\mathcal{B}}\mu)\Vert \bar{x}(k)-x^{\star}\Vert_{2}.\notag
\end{eqnarray}
\end{lemma}

The first and third inequalities in Lemma \ref{Found1} follow from Assumption \ref{LipschitzF} and the fact that $\Vert \bar{g}(k)\Vert_{2}=\Vert \bm{1}_{n}^{\text{T}}\nabla F(\bm{1}_{n}\bar{x}(k))-\bm{1}_{n}^{\text{T}}\nabla F(\bm{1}_{n}x^{\star})\Vert_{2}$,
while the second inequality can be obtained directly by applying (\ref{Important111}). 
The last statement can be verified by following the similar line of Lemma 10 in \cite{QuG2018Harnssing}. 
Now, we begin to prove Lemma \ref{proposition1} by establishing the upper bounds of $\Vert \bar{x}(k+1)-x^{\star}\Vert_{2}$,  $\Vert x(k+1)-\bm{1}_{n}\bar{x}(k+1)\Vert_{\mathcal{A}}$ and $\Vert z(k+1)-v\bar{z}(k+1)\Vert_{\mathcal{B}}$, respectively. 
%

(i) 
In view of (\ref{barxx}), we have 
\begin{eqnarray}
\bar{x}(k+1)
&\hspace*{-0.5em}=\hspace*{-0.5em}&
\bar{x}(k)-\eta \pi_{\mathcal{A}}^{\text{T}}\pi_{\mathcal{B}}\bar{z}(k)-\eta \pi_{\mathcal{A}}^{\text{T}}
(z(k)-\pi_{\mathcal{B}}\bar{z}(k))  \notag\\
&\hspace*{-0.5em}=\hspace*{-0.5em}&
\bar{x}(k)
-\eta \pi_{\mathcal{A}}^{\text{T}}\pi_{\mathcal{B}}\bar{g}(k)
-\eta\pi_{\mathcal{A}}^{\text{T}}\pi_{\mathcal{B}}(\bar{z}(k)-g(k))
  \notag\\
&&-\eta\pi_{\mathcal{A}}^{\text{T}}\pi_{\mathcal{B}}(g(k)-\bar{g}(k))
-\eta\pi_{\mathcal{A}}^{\text{T}}(z(k)-\pi_{\mathcal{B}}\bar{z}(k))  \notag\\
&\hspace*{-0.5em}=\hspace*{-0.5em}&
\bar{x}(k)
-\eta \pi_{\mathcal{A}}^{\text{T}}\pi_{\mathcal{B}}\bar{g}(k)
-\eta \pi_{\mathcal{A}}^{\text{T}}\pi_{\mathcal{B}}(g(k)-\bar{g}(k))  \notag\\
&&-\eta \pi_{\mathcal{A}}^{\text{T}}(z(k)-\pi_{\mathcal{B}}\bar{z}(k))
-\eta  \pi_{\mathcal{A}}^{\text{T}}\pi_{\mathcal{B}}\bm{1}_{n}^{\text{T}}\epsilon_{y}(k-1). \notag
\end{eqnarray}
Therefore, by invoking Lemma \ref{Normequ2} and Lemma \ref{Found1}, we further obtain 
\begin{eqnarray}\label{AuxL1}
\hspace*{-1em}&&\Vert \bar{x}(k+1)-x^{\star}\Vert_{2}\notag\\
\hspace*{-1em}&&\leq
(1-\eta \pi_{\mathcal{A}}^{\text{T}}\pi_{\mathcal{B}}\mu)\Vert \bar{x}(k)-x^{\star}\Vert_{2}
+\eta \pi_{\mathcal{A}}^{\text{T}}\pi_{\mathcal{B}}\Vert g(k)-\bar{g}(k))\Vert_{2}  \notag\\
\hspace*{-1em}&&~~~+\eta\Vert \pi_{\mathcal{A}}^{\text{T}}(z(k)-\pi_{\mathcal{B}}\bar{z}(k))\Vert_{2}
+\eta \pi_{\mathcal{A}}^{\text{T}}\pi_{\mathcal{B}}\left\Vert \bm{1}_{n}^{\text{T}}\epsilon_{y}(k-1)\right\Vert_{2} \notag\\
\hspace*{-1em}&&\leq
(1-\eta \pi_{\mathcal{A}}^{\text{T}}\pi_{\mathcal{B}}\mu)\Vert \bar{x}(k)-x^{\star}\Vert_{2}
+\sqrt{n}\eta \pi_{\mathcal{A}}^{\text{T}}\pi_{\mathcal{B}}\Vert x(k)-\bm{1}_{n}\bar{x}(k)\Vert_{\mathcal{A}}  \notag\\
\hspace*{-1em}&&~~~+\eta \Vert z(k)-\pi_{\mathcal{B}}\bar{z}(k)\Vert_{\mathcal{B}}
+\eta \pi_{\mathcal{A}}^{\text{T}}\pi_{\mathcal{B}}\Vert\bm{1}_{n}^{\text{T}}\epsilon_{y}(k-1)\Vert_{2}.
\end{eqnarray}
where the fact that $\Vert \pi_{\mathcal{A}}\Vert_{2}\leq 1$ has been used to obtain the last inequality. 
(ii)
From (\ref{AlgRes_1A}) and (\ref{barxx}), along with $\mathcal{A}_{\alpha}\bm{1}_{n}=\bm{1}_{n}$, we obtain 
\begin{eqnarray}
&&x(k+1)-\bm{1}_{n}\bar{x}(k+1)  \notag\\
&&=
\mathcal{A}_{\alpha}(x(k)-\bm{1}_{n}\bar{x}(k))-\eta (I_{n}-\bm{1}_{n}\pi_{\mathcal{A}}^{\text{T}})z(k) \notag\\
&&~~~+\alpha(\mathcal{A}-I_{n})\sigma_{x}(k)  \notag\\
&&=
(\mathcal{A}_{\alpha}-\bm{1}_{n}\pi_{\mathcal{A}}^{\text{T}})(x(k)-\bm{1}_{n}\bar{x}(k))
-\eta (I_{n}-\bm{1}_{n}\pi_{\mathcal{A}}^{\text{T}})z(k) \notag\\
&&~~~+\alpha(\mathcal{A}-I_{n})\sigma_{x}(k)  \notag\\
&&=
(\mathcal{A}_{\alpha}-\bm{1}_{n}\pi_{\mathcal{A}}^{\text{T}})(x(k)-\bm{1}_{n}\bar{x}(k))
-\eta (I_{n}-\bm{1}_{n}\pi_{\mathcal{A}}^{\text{T}})\pi_{\mathcal{B}}\bar{z}(k) \notag\\
&&~~~-\eta (I_{n}-\bm{1}_{n}\pi_{\mathcal{A}}^{\text{T}})(z(k)
-\pi_{\mathcal{B}}\bar{z}(k))
+\alpha(\mathcal{A}-I_{n})\sigma_{x}(k),  \notag
\end{eqnarray}
where  the fact that $\bm{1}_{n}\pi_{\mathcal{A}}^{\text{T}}(x(k)-\bm{1}_{n}\bar{x}(k))=\bm{0}_{n}$ has been exploited to obtain the second equality.
By employing Lemma \ref{Normequ2}, we obtain 
\begin{eqnarray}\label{Imp112}
&&\Vert x(k+1)-\bm{1}_{n}\bar{x}(k+1)\Vert_{\mathcal{A}}   \notag\\
&&\leq
\sigma_{\mathcal{A}}\Vert x(k)-\bm{1}_{n}\bar{x}(k)\Vert_{\mathcal{A}}
+\alpha\Vert(\mathcal{A}-I_{n})\sigma_{x}(k)\Vert_{\mathcal{A}}           \notag\\
&&~~~+\eta \kappa_{1}
\Vert z(k)-\pi_{\mathcal{B}}\bar{z}(k) \Vert_{\mathcal{A}}
+\eta\kappa_{1}\kappa_{2} \delta_{\mathcal{A}2}\Vert \bar{z}(k) \Vert_{2}. 
\end{eqnarray}
Now, it remains to establish an upper bound for $\Vert \bar{z}(k) \Vert_{2}$.
Note that
\begin{eqnarray}\label{Imp111}
\Vert \bar{z}(k) \Vert_{2}
&\hspace*{-0.5em}\leq\hspace*{-0.5em}&
\Vert \bar{z}(k)-g(k) \Vert_{2} +\Vert g(k)-\bar{g}(k) \Vert_{2}+\Vert \bar{g}(k) \Vert_{2}  \notag\\
&\hspace*{-0.5em}\leq\hspace*{-0.5em}&
 \Vert\bm{1}_{n}^{\text{T}}\epsilon_{y}(k-1)\Vert_{2}
+\sqrt{n}L\Vert x(k)-\bm{1}_{n}\bar{x}(k)\Vert_{\mathcal{A}}  \notag\\
&&+nL\Vert \bar{x}(k)-x^{\star}\Vert_{2}.
\end{eqnarray}
By substituting (\ref{Imp111}) into (\ref{Imp112}), 
we obtain 
\begin{eqnarray}\label{IMP22}
\hspace*{-1em}&&\Vert x(k+1)-\bm{1}_{n}\bar{x}(k+1)\Vert_{\mathcal{A}}   \notag\\
\hspace*{-1em}&&\leq
\left(\sigma_{\mathcal{A}} +\sqrt{n}\eta L \kappa_{1}\kappa_{2}\delta_{\mathcal{A}2}\right)
\Vert x(k)-\bm{1}_{n}\bar{x}(k)\Vert_{\mathcal{A}}
         \notag\\
\hspace*{-1em}&&~~~+\eta\delta_{\mathcal{AB}} \kappa_{1}
\Vert z(k)-\pi_{\mathcal{B}}\bar{z}(k) \Vert_{\mathcal{B}}
+n\eta L \kappa_{1}\kappa_{2}\delta_{\mathcal{A}2}\Vert \bar{x}(k)-x^{\star} \Vert_{2}   \notag\\
\hspace*{-1em}&&~~~+\eta \kappa_{1}\kappa_{2}\delta_{\mathcal{A}2} \Vert\bm{1}_{n}^{\text{T}}\epsilon_{y}(k-1)\Vert_{2}   
+\alpha\Vert(\mathcal{A}-I_{n})\sigma_{x}(k)\Vert_{\mathcal{A}}. 
\end{eqnarray}


(iii)
In light of relations (\ref{AlgRes_2B}) and (\ref{yyyya}), we have
\begin{eqnarray}
&&z(k+1)-\pi_{\mathcal{B}}\bar{z}(k+1)      \notag\\
&&=
(\mathcal{B}_{\beta}-\pi_{\mathcal{B}}\bm{1}_{n}^{\text{T}})\left(z(k)-\pi_{\mathcal{B}}\bar{z}(k)\right)
 \notag\\
&&~~~+(I_{n}-\pi_{\mathcal{B}}\bm{1}_{n}^{\text{T}})(\nabla F(x(k+1))-\nabla F(x(k)))  \notag\\
&&~~~+(I_{n}-\pi_{\mathcal{B}}\bm{1}_{n}^{\text{T}})(\epsilon_{y}(k) -\epsilon_{y}(k-1)), \notag
\end{eqnarray}
where the equality follows from the definition of $\mathcal{B}_{\beta}$ and the fact that $\pi_{\mathcal{B}}\bm{1}_{n}^{\text{T}}\pi_{\mathcal{B}}=\pi_{\mathcal{B}}$.
Hence, 
we obtain 
\begin{eqnarray}\label{ddddf}
&&\Vert z(k+1)-\pi_{\mathcal{B}}\bar{z}(k+1)\Vert_{\mathcal{B}} \notag\\
&&\leq
\delta_{\mathcal{B}2}\Vert I_{n}-\pi_{\mathcal{B}}\bm{1}_{n}^{\text{T}}\Vert_{\mathcal{B}} \Vert \nabla F(x(k+1))-\nabla F(x(k))\Vert_{2} \notag \\
&&~~~+\Vert I_{n}-\pi_{\mathcal{B}}\bm{1}_{n}^{\text{T}}\Vert_{\mathcal{B}}\Vert \epsilon_{y}(k)-\epsilon_{y}(k-1) \Vert_{\mathcal{B}}\notag\\
&&~~~+\sigma_{\mathcal{B}}\Vert z(k)-\pi_{\mathcal{B}}\bar{z}(k)\Vert_{\mathcal{B}}, 
\end{eqnarray}
%
where Lemma \ref{Normequ2} has been utilized to obtain the above inequality.
Now, it remains to bound $\Vert \nabla F(x(k+1))-\nabla F(x(k))\Vert_{2}$.
Note that
\begin{eqnarray}\label{dsss}
&&\Vert \nabla F(x(k+1))-\nabla F(x(k))\Vert_{2} \notag\\
&&\leq
L\Vert x(k+1)-x(k)\Vert_{2} \notag \\
&&=
L\Vert \mathcal{A}_{\alpha}x(k)-x(k)+\alpha(\mathcal{A}-I_{n})\sigma_{x}(k)-\eta z(k)\Vert_{2}\notag \\
&&\leq
L\Vert \mathcal{A}_{\alpha}-I_{n}\Vert_{2}\Vert x(k)-\bm{1}_{n}\bar{x}(k)\Vert_{2}
+\eta L\Vert z(k)-\pi_{\mathcal{B}}\bar{z}(k)  \Vert_{2}
\notag \\
&&~~~+\eta L \Vert  \bar{z}(k)\Vert_{2}
+\alpha L \Vert(\mathcal{A}-I_{n})\sigma_{x}(k)\Vert_{2}. 
\end{eqnarray} 
where the fact that $\Vert \pi_{\mathcal{B}}\Vert_{2}\leq 1$ has been used to obtain the last inequality. 
Then, by substituting (\ref{Imp111}) and (\ref{dsss}) into (\ref{ddddf}), we can obtain 
\begin{eqnarray}\label{IMPwww}
&&\Vert z(k+1)-\pi_{\mathcal{B}}\bar{z}(k+1)\Vert_{\mathcal{B}} \notag\\
&&\leq
\delta_{\mathcal{B}2}\kappa_{3}\left(L\kappa_{4}+\sqrt{n}\eta L^{2}\right)\Vert x(k)-\bm{1}_{n}\bar{x}(k)\Vert_{\mathcal{A}}\notag\\
&&~~~+\delta_{\mathcal{B}2}\alpha L\kappa_{3}\Vert (\mathcal{A}-I_{n})\sigma_{x}(k)\Vert_{2}  
+\delta_{\mathcal{B}2}\eta L\kappa_{3}\Vert\bm{1}_{n}^{\text{T}}\epsilon_{y}(k-1)\Vert_{2}\notag\\
&&~~~+\kappa_{3}\Vert \epsilon_{y}(k)- \epsilon_{y}(k-1) \Vert_{\mathcal{B}}
+\delta_{\mathcal{B}2}n\eta L^{2}\kappa_{3}\Vert \bar{x}(k)-x^{\star}\Vert_{2}  \notag\\
&&~~~+(\sigma_{\mathcal{B}}+\delta_{\mathcal{B}2}\eta L\kappa_{3})
\Vert z(k)-\pi_{\mathcal{B}}\bar{z}(k)\Vert_{\mathcal{B}}. 
\end{eqnarray}

Combining (\ref{AuxL1}), (\ref{IMP22}) and (\ref{IMPwww}), we can obtain (\ref{System}) with
\begin{align}\label{GmatrixGG}
	G=
	\left[
	\begin{array}{ccc}
		1-\eta \pi_{\mathcal{A}}^{\text{T}}\pi_{\mathcal{B}}\mu & \sqrt{n}\eta \pi_{\mathcal{A}}^{\text{T}}\pi_{\mathcal{B}} & \eta \\
		n\eta L \kappa_{1}\kappa_{2}\delta_{\mathcal{A}2} & \sigma_{\mathcal{A}} +\sqrt{n}\eta L  \kappa_{1}\kappa_{2}\delta_{\mathcal{A}2} & \eta \kappa_{1}\delta_{\mathcal{AB}} \\
		n\eta L^{2}\kappa_{3}\delta_{\mathcal{B}2} & \kappa_{3}
		(L\kappa_{4}+\sqrt{n}\eta L^{2})\delta_{\mathcal{B}2} & \sigma_{\mathcal{B}}+\eta L\kappa_{3}\delta_{\mathcal{B}2} \\
	\end{array}
	\right]
\end{align}
and $\varsigma(k) = (\varsigma_1(k), \varsigma_2(k), \varsigma_3(k))^{\text{T}}$ given by 
\begin{eqnarray}\label{varsigma123}
	\varsigma_{1}(k)
	&\hspace*{-0.5em}=\hspace*{-0.5em}&
	\eta \pi_{\mathcal{A}}^{\text{T}}\pi_{\mathcal{B}}\Vert\bm{1}_{n}^{\text{T}}\epsilon_{y}(k-1)\Vert_{2}, \notag\\
	\varsigma_{2}(k)
	&\hspace*{-0.5em}=\hspace*{-0.5em}&
	\eta \kappa_{1}\kappa_{2}\delta_{\mathcal{A}2} \Vert\bm{1}_{n}^{\text{T}}\epsilon_{y}(k-1)\Vert_{2}
	+\alpha\Vert(\mathcal{A}-I_{n})\sigma_{x}(k)\Vert_{\mathcal{A}},   \notag\\
	\varsigma_{3}(k)
	&\hspace*{-0.5em}=\hspace*{-0.5em}&
	\kappa_{3}\Vert \epsilon_{y}(k)-\epsilon_{y}(k-1) \Vert_{\mathcal{B}}
	+\eta L\kappa_{3}\delta_{\mathcal{B}2}\Vert\bm{1}_{n}^{\text{T}}\epsilon_{y}(k-1)\Vert_{2} \notag\\
		&&+\alpha L\kappa_{3}\delta_{\mathcal{B}2}\Vert (\mathcal{A}-I_{n})\sigma_{x}(k)\Vert_{2},  
\end{eqnarray}
which completes the proof. 
\hfill $\blacksquare$


\vspace{-1ex}

\subsection{\textit{Proof of Lemma \ref{StepCond}.}}
To achieve this goal, 
we need to provide a sufficient condition 
under which $G_{ii} < 1$ and $\det(I - G) > 0$ can be guaranteed \cite{PuS2018PushJou}. 
 We first ensure that $G_{ii}<1$ hold for $i=1,2,3$.
Clearly, 
if we set $\eta\leq \frac{1}{(\mu+L)\pi_{\mathcal{A}}^{\text{T}}\pi_{\mathcal{B}}}$, then $0<G_{11}<1$. 
We can also verify that if 
$\eta \leq \min \{\frac{1-\sigma_{\mathcal{A}}}{2\sqrt{n}L\kappa_1\kappa_2 \delta_{\mathcal{A}2}}, \frac{1-\sigma_{\mathcal{B}}}{2L\kappa_{3}\delta_{\mathcal{B}2}}\}$,
then $G_{22}<1$ and $G_{33}<1$ both hold.
Now we turn our attention to $\det(I_{3}- G)$.
Note that 
\begin{eqnarray}
&&	\det(I_{3} - G) \notag\\
	&&=
	\eta \pi_{\mathcal{A}}^{\text{T}}\pi_{\mathcal{B}} \mu   [1 - (\sigma_{\mathcal{A}}+\sqrt{n} \eta L \kappa_{1} \kappa_{2}\delta_{\mathcal{A}2})]  [1 - (\sigma_{\mathcal{B}}+ \eta L \kappa_{3}\delta_{\mathcal{B} 2})] \notag \\
		&&~~~ - \eta \pi_{\mathcal{A}}^{\text{T}}\pi_{\mathcal{B}} \mu  (\eta \kappa_{1}\delta_{\mathcal{AB}} )  [ \kappa_{3}(L \kappa_{4}+\sqrt{n} \eta L^{2})\delta_{\mathcal{B} 2}] \notag\\
			&&~~~  - \sqrt{n}\eta \pi_{\mathcal{A}}^{\text{T}}\pi_{\mathcal{B}}  (n \eta L \kappa_{1} \kappa_{2}\delta_{\mathcal{A}2})  [1 - (\sigma_{\mathcal{B}}+ \eta L \kappa_{3}\delta_{\mathcal{B} 2})] \notag\\
	&&  ~~~ - \sqrt{n} \eta \pi_{\mathcal{A}}^{\text{T}}\pi_{\mathcal{B}}  (\eta\kappa_{1} \delta_{\mathcal{AB}} )  ( n \eta L^{2} \kappa_{3}\delta_{\mathcal{B} 2} )   \notag\\
	&&~~~ - \eta   (n \eta L \kappa_{1} \kappa_{2}\delta_{\mathcal{A}2})  [ \kappa_{3}(L \kappa_{4}+\sqrt{n} \eta L^{2} )\delta_{\mathcal{B} 2}] \notag\\
	&&~~~  - \eta   [1 - (\sigma_{\mathcal{A}}+\sqrt{n} \eta L \kappa_{1} \kappa_{2}\delta_{\mathcal{A}2})]  ( n \eta L^{2} \kappa_{3}\delta_{\mathcal{B} 2}). \notag
\end{eqnarray}
In light of $\eta \leq \min \{\frac{1-\sigma_{\mathcal{A}}}{2\sqrt{n}L\kappa_1\kappa_2\delta_{\mathcal{A}2}}, \frac{1-\sigma_{\mathcal{B}}}{2L\kappa_{3}\delta_{\mathcal{B}2}}\}$, we have
$\frac{1 - \sigma_{\mathcal{A}}}{2}\leq 1 - (\sigma_{\mathcal{A}}+\sqrt{n} \eta L \kappa_{1} \kappa_{2} \delta_{\mathcal{A}2})
\leq 1 - \sigma_{\mathcal{A}}$
and $\frac{1 - \sigma_{\mathcal{B}}}{2}\leq
 1 - (\sigma_{\mathcal{B}}+\delta_{\mathcal{B} 2} \eta L \kappa_{3})\leq 1 - \sigma_{\mathcal{B}}$.
Therefore, a sufficient condition for $\det(I - G) > 0$ is
\begin{eqnarray}\label{detgeq0}
	&& \frac{1}{4}\eta \pi_{\mathcal{A}}^{\text{T}}\pi_{\mathcal{B}} \mu  (1 - \sigma_{\mathcal{A}}) (1 - \sigma_{\mathcal{B}})- \eta  (1 - \sigma_{\mathcal{A}})  ( n \eta L^{2} \kappa_{3} \delta_{\mathcal{B} 2})\notag\\
	&& - \eta \pi_{\mathcal{A}}^{\text{T}}\pi_{\mathcal{B}}  \mu (\eta  \kappa_{1}\delta_{\mathcal{A} \mathcal{B}})  [\kappa_{3}(L \kappa_{4}+\sqrt{n} \eta L^{2})\delta_{\mathcal{B} 2} ]\notag\\
	&&	- \sqrt{n} \eta \pi_{\mathcal{A}}^{\text{T}}\pi_{\mathcal{B}} (n \eta L \kappa_{1} \kappa_{2}\delta_{\mathcal{A}2})  (1 - \sigma_{\mathcal{B}}) \notag\\
	&&	- \sqrt{n} \eta \pi_{\mathcal{A}}^{\text{T}}\pi_{\mathcal{B}}   (\eta  \kappa_{1}\delta_{\mathcal{A} \mathcal{B}})  ( n \eta L^{2} \kappa_{3}\delta_{\mathcal{B} 2}) \notag\\
	&& - \eta  (n \eta L \kappa_{1} \kappa_{2}\delta_{\mathcal{A}2})  [ \kappa_{3}(L \kappa_{4}+\sqrt{n} \eta L^{2} )\delta_{\mathcal{B} 2}] > 0. 
\end{eqnarray}
%
%
Now, the inequality (\ref{detgeq0}) can be rewritten as $\Gamma_1\eta^2 + \Gamma_2\eta -\Gamma_3 < 0$ with 
$\Gamma_{i}$, $i=1,2,3$, given by 
\begin{eqnarray}\label{GXi1}
		\Gamma_{1}
	&\hspace*{-0.5em}=\hspace*{-0.5em}&
	\sqrt{n}\kappa_{1}\kappa_{3} L^{2}\delta_{\mathcal{B}2}[(n+\mu)\pi_{\mathcal{A}}^{\text{T}} \pi_{\mathcal{B}}\delta_{\mathcal{AB}}+n\kappa_{2}L\delta_{\mathcal{A}2}],\notag \\
\Gamma_{2}
&\hspace*{-0.5em}=\hspace*{-0.5em}&
\kappa_{1}L\pi_{\mathcal{A}}^{\text{T}} \pi_{\mathcal{B}}[n^{\frac{3}{2}}\kappa_{2}\delta_{\mathcal{A}2}(1-\sigma_{\mathcal{B}})+\mu\kappa_{3}\kappa_{4} \delta_{\mathcal{AB}}\delta_{\mathcal{B}2}]\notag\\
&&+n\kappa_{3}L^{2}\delta_{\mathcal{B}2} [1-\sigma_{\mathcal{A}}+\kappa_{1}\kappa_{2}\kappa_{4}\delta_{\mathcal{A}2}],\notag \\
		\Gamma_{3}
		&\hspace*{-0.5em}=\hspace*{-0.5em}&
		\frac{1}{4}\mu \pi_{\mathcal{A}}^{\text{T}} \pi_{\mathcal{B}}(1-\sigma_{\mathcal{A}})(1-\sigma_{\mathcal{B}}), 
	\end{eqnarray}
%
Therefore, it can be derived that
$\eta \leq  \frac{2\Gamma_3}{\Gamma_2+\sqrt{\Gamma_2^2+4\Gamma_1\Gamma_3}}$, which completes the proof.
\hfill $\blacksquare$

%
%
%

%
%

\vspace{-1ex}

\subsection{\textit{Proof of Theorem \ref{Unsat1}.}}

To ensure that the finite-level quantizers never saturate, 
	the scaled ``innovation'' $\frac{1}{h(k)}(x_{j}(k)-\hat{x}_{j}(k-1))$ and $\frac{1}{h(k)}(y_{j}(k)-\hat{y}_{j}(k-1))$ must lie in a bounded region.
	To achieve the goal, we first establish the upper bounds for $\frac{1}{h(k)}\Vert x_{i}(k)-\hat{x}_{i}(k-1)\Vert_{\infty}$ and $ \frac{1}{h(k)}\Vert y_{i}(k)-\hat{y}_{i}(k-1)\Vert_{\infty}$, respectively. 
	Then, the obtained upper bounds lead us to propose an update rule of the quantization levels, under which we prove the unsaturation of quantizers by mathematical induction. 
	Finaly, we show that (\ref{QuanLevelConAlterQQQ}) suffices for the given update rule. \\
%
%
%
%
\textbf{Step 1: Bound $\Vert x_{i}(k)-\hat{x}_{i}(k-1)\Vert_{\infty}$ and
$ \Vert y_{i}(k)-\hat{y}_{i}(k-1)\Vert_{\infty}$.}

Let 
$e_{x_{i}}(k)\triangleq Q_{K_{x}}(\frac{x_{i}(k)-\hat{x}_{i}(k-1)}{h(k)})-\frac{x_{i}(k)-\hat{x}_{i}(k-1)}{h(k)}$
and $e_{y_{i}}(k)\triangleq Q_{K_{y}}(\frac{y_{i}(k)-\hat{y}_{i}(k-1)}{h(k)})-\frac{y_{i}(k)-\hat{y}_{i}(k-1)}{h(k)}$. 
Recalling (\ref{rrrrjk}) and (\ref{hatttttxjk}), we can obtain $\hat{x}_{j}(k)=x_{j}(k)+h(k)e_{x_{j}}(k)$. Then  
%
\begin{eqnarray}\label{BoundImp1}
&& \Vert x_{i}(k)-\hat{x}_{i}(k-1)\Vert_{\infty} \notag\\
&&\leq
  \Vert x_{i}(k)-x_{i}(k-1)\Vert_{\infty}
+h(k-1)\left\Vert e_{x_i}(k-1)\right\Vert_{\infty}. 
\end{eqnarray}
The first term on the right side of (\ref{BoundImp1}) can be further calculated as
%
%
%
 \begin{eqnarray}\label{ThreeTerm}
&&\Vert x_i(k) - x_i(k-1)\Vert_{\infty}  
\leq
 \alpha\sum_{j=1}^{n}\Vert x_{i}(k-1)-x_{j}(k-1)\Vert_{\infty}   \notag\\
&& +\eta\Vert z_{i}(k-1)\Vert_{\infty} 
+ \alpha\sum_{j=1}^{n}\Vert \sigma_{x_{i}}(k-1)-\sigma_{x_{j}}(k-1)\Vert_{\infty}, 
 \end{eqnarray}
where the inequality follows from (\ref{AlgRes_1A}), the definition of $\mathcal{A}_{\alpha}$, and the row stochasticity of $\mathcal{A}$.  
In the following, we will establish the upper bounds for the three 
terms on the right side of (\ref{ThreeTerm}), 
respectively. 

For the first term, it can be calculated as follows
\begin{eqnarray}\label{}
\sum_{j=1}^{n}\Vert x_{i}(k-1)-x_{j}(k-1)\Vert_{\infty}
 &\hspace*{-0.5em}\leq\hspace*{-0.5em}&
 \sqrt{2}(n+\frac{1}{2})\Theta_{2}(k-1), 
\end{eqnarray}
where 
the Jensen's inequality and the facts that
$\Vert x_{i}(k-1)-x_{j}(k-1)\Vert^{2}_{\infty}\leq 2\Vert x_{i}(k-1)-\bar{x}(k-1)\Vert^{2}_{\infty}
+2\Vert \bar{x}(k-1)-x_{j}(k-1)\Vert^{2}_{\infty}$ 
and
$\Vert x_{i}(k-1)-\bar{x}(k-1)\Vert^{2}_{\infty}\leq\Vert x_{i}(k-1)-\bar{x}(k-1)\Vert^{2}_{2}$ 
have been exploited to obtain the above inequality.

For the second term on the right side of (\ref{ThreeTerm}), we have 
\begin{eqnarray}\label{}
&&\Vert z_{i}(k-1)\Vert_{\infty}  \notag\\
&&\leq
 \left\Vert z(k-1) - \pi_{\mathcal{B}}\bar{z}(k-1)\right\Vert_{\infty} 
 +\Vert\pi_{\mathcal{B}}\Vert_{\infty}\left\Vert \bar{g}(k-1) \right\Vert_{\infty} \notag\\
&&~~~+\Vert\pi_{\mathcal{B}}\Vert_{\infty} \left\Vert \bar{z}(k-1) - g(k-1) \right\Vert_{\infty}   \notag\\
&&~~~+\Vert\pi_{\mathcal{B}}\Vert_{\infty}\left\Vert g(k-1) - \bar{g}(k-1) \right\Vert_{\infty}  \notag\\
&&\leq
\Theta_{3}(k-1)+\sqrt{n}L\Theta_{2}(k-1)+nL\Theta_{1}(k-1)  \notag\\
&&~~~+n\beta h(k-2)\max_{i\in\mathcal{V}}\Vert e_{y_{i}}(k-2)\Vert_{\infty}, 
\end{eqnarray}
where  
Lemma \ref{Normequ2} and Lemma \ref{Found1} have been employed to obtain the above inequality.

It only remains to bound the last term in (\ref{ThreeTerm}).  
By using the fact that $\hat{x}_{j}(k)=x_{j}(k)+h(k)e_{x_{j}}(k)$ again, we obtain  
$\sum_{j=1}^{n}\Vert \sigma_{x_{i}}(k-1)-\sigma_{x_{j}}(k-1)\Vert_{\infty}\leq 2 n h(k-1)\max_{i\in\mathcal{V}}\Vert e_{x_{i}}(k-1) \Vert_{\infty}$. 
%

Define 
\begin{eqnarray}\label{varpi123}
	\varphi_{1}
	&\hspace*{-0.5em}\triangleq\hspace*{-0.5em}&	
	\max\left\{\sqrt{2}(n+\frac{1}{2})\alpha+\eta\sqrt{n}L, \eta, \eta nL\right\},\notag \\
	\varphi_{2}
	&\hspace*{-0.5em}\triangleq\hspace*{-0.5em}&
	\max\left\{1, \sqrt{n}L,  nL\right\}.  
\end{eqnarray}
Combining the above inequalities, we can obtain 
\begin{eqnarray}\label{IMPxmihatx}
&&\Vert x_{i}(k)-\hat{x}_{i}(k-1)\Vert_{\infty} 
\leq
\sqrt{3}\varphi_{1}\Vert \Theta(k-1) \Vert_{2}  \notag\\
&&~~~~~~~~~~~~~~+(2\alpha n+1)h(k-1)\max_{i\in\mathcal{V}}\Vert e_{x_{i}}(k-1)\Vert_{\infty}\notag\\
&&~~~~~~~~~~~~~~+n \eta \beta h(k-2)\max_{i\in\mathcal{V}}\Vert e_{y_{i}}(k-2)\Vert_{\infty}.
\end{eqnarray}
From (\ref{IMPxmihatx}), we can observe that if the quantizers never saturate, then $x_{j}(k)-\hat{x}_{j}(k-1)$ will decay to zero 
	at the speed of the 
	same order of $h(k)$ 
	since $h(k-1)=\frac{C}{\xi}\xi^{k}$. 
Following the similar line above, we can further obtain 
%
%
%
%
%
\begin{eqnarray}\label{}
 &&\Vert y_{i}(k)-\hat{y}_{i}(k-1)\Vert_{\infty}
\leq
\sqrt{3}\varphi_{2}\Vert \Theta(k) \Vert_{2}   \notag\\
&&~~~~~~~~~~~~+(n \beta +1) h(k-1)\max_{i\in\mathcal{V}}\Vert e_{y_{i}}(k-1)\Vert_{\infty}. \notag
\end{eqnarray}
\textbf{Step 2: Demonstrate the unsaturation.}

	In this part, we first consider the following update rule of the quantization levels instead 
	\begin{eqnarray}\label{QuanLevelConAlter}
		K_x(0) &\hspace*{-0.8em}\geq\hspace*{-0.8em}& \frac{v_{1}}{C} - \frac{1}{2}, ~
		K_y(0) \geq \frac{v_{2}}{C} - \frac{1}{2} \notag\\
		K_x(1) &\hspace*{-0.8em}\geq\hspace*{-0.8em}& \frac{\sqrt{3}\varphi_{1}\Vert\Theta(0)\Vert_{2}}{C \xi} + \frac{2\alpha n+1}{2\xi}  - \frac{1}{2} \notag \\
		K_x(k) &\hspace*{-0.8em}\geq\hspace*{-0.8em}& \frac{\sqrt{3}\varphi_{1} \tau \Vert\Theta(0)\Vert_{2}}{C \xi} \Upsilon_1(k) + \frac{2\alpha n+1}{2\xi} + \frac{ n\eta\beta }{2  \xi^2} - \frac{1}{2}, k \geq 2 \notag\\
		K_y(k) &\hspace*{-0.8em}\geq\hspace*{-0.8em}& \frac{\sqrt{3}\varphi_{2}\tau \Vert\Theta(0)\Vert_{2}}{C} \Upsilon_2(k) + \frac{n \beta +1 }{2\xi} - \frac{1}{2}, ~k \geq 1
	\end{eqnarray}
	where 
	$\Upsilon_1(k)=(\frac{\hat{\rho}}{\xi})^{k-1} + \frac{\tilde{\varsigma}}{\xi \Vert\Theta(0)\Vert_{2}}\sum_{l=0}^{k-3}(\frac{\hat{\rho}}{\xi})^{k-2-l} +\frac{\tilde{\varsigma}}{\xi\tau \Vert\Theta(0)\Vert_{2}}$
	and 
	$\Upsilon_2(k)=(\frac{\hat{\rho}}{\xi})^{k} + \frac{\tilde{\varsigma}}{\xi \Vert\Theta(0)\Vert_{2}}\sum_{l=0}^{k-2}(\frac{\hat{\rho}}{\xi})^{k-1-l} +\frac{\tilde{\varsigma}}{\xi\tau \Vert\Theta(0)\Vert_{2}}$. 

Now, we show the unsaturation of the quantizers under the rule (\ref{QuanLevelConAlter}) by mathematical induction.  
Considering the case $k=0$, we have
$ \frac{\Vert x_{i}(k)-\hat{x}_{i}(k-1)\Vert_{\infty}}{h(k)}
\leq \frac{\Vert x_{i}(0)\Vert_{\infty}}{C}
\leq K_{x}(0)+\frac{1}{2}$
and
$\frac{\Vert y_{i}(k)-\hat{y}_{i}(k-1)\Vert_{\infty}}{h(k)} \leq \frac{\Vert y_{i}(0)\Vert_{\infty}}{C}
\leq K_{y}(0)+\frac{1}{2}$,
which indicates that
the quantizers are not saturated
for $k=0$. Therefore, $\max_{i \in \mathcal{V}} \Vert e_{x_{i}}(0)\Vert_{\infty} \leq \frac{1}{2}$ and $\max_{i \in \mathcal{V}} \Vert e_{y_{i}}(0) \Vert_{\infty}\leq \frac{1}{2}$ both hold, which further can be exploited to calculate the upper bounds of
$\varsigma_{i}(0)$ via (\ref{varsigma123}), denoted by $\bar{\varsigma}_{i}(0)$, for $i=1,2,3$.
Define $\hat{\varsigma}(0) \triangleq\Vert \bar{\varsigma}(0)\Vert_{2}$ with $\bar{\varsigma}(0)=(\bar{\varsigma}_{1}(0),\bar{\varsigma}_{2}(0),\bar{\varsigma}_{3}(0))^{\text{T}}$.
Recalling (\ref{System}), we can obtain
$\Vert \Theta(1)\Vert_{2}\leq \tau \hat{\rho}\Vert\Theta(0)\Vert_{2}+ \hat{\varsigma}(0)$. 

Now, considering the case $k=1$. From (\ref{IMPxmihatx}), we can obtain
$\frac{\Vert x_{i}(1)-\hat{x}_{i}(0)\Vert_{\infty}}{h(1)} \leq
\frac{\sqrt{3}\varphi_{1}}{C \xi}\Vert \Theta(0) \Vert_{2}
+\frac{2\alpha n+1}{\xi}\max_{i\in\mathcal{V}}\Vert e_{x_{i}}(0)\Vert_{\infty}
\leq
\frac{\sqrt{3}\varphi_{1}}{C \xi}\Vert \Theta(0) \Vert_{2} + \frac{2\alpha n+1}{2\xi}
\leq
 K_x(1) + \frac{1}{2}$.
 Similarly, it can be easily verified that $\frac{\Vert y_{i}(1)-\hat{y}_{i}(0)\Vert_{\infty}}{h(1)}\leq K_y(1) + \frac{1}{2}$.
These two inequalities imply that the quantizers are not saturated at $k=1$ as well. 
Then, 
we have 
$\max_{i \in \mathcal{V}} \Vert e_{x_{i}}(\nu) \Vert_{\infty} \leq \frac{1}{2}$ and $\max_{i \in \mathcal{V}} \Vert e_{y_{i}}(\nu)\Vert_{\infty} \leq \frac{1}{2}$ for $\nu \in \{0,1\}$,
which further can be utilized to compute $\hat{\varsigma}(1)$.
Hence, we can obtain  $\Vert\Theta(2)\Vert_2
\leq
\tau \Vert \Theta(0) \Vert_{2} \hat{\rho}^2 + \tau\hat{\rho}\hat{\varsigma}(0) + \hat{\varsigma}(1)$.

From the above observations, it can be seen that our basic idea is to exploit the non-saturation property at each step,   
i.e., $\max_{i \in \mathcal{V}} \Vert e_{x_i}(\nu)\Vert_{\infty} \leq \frac{1}{2}$ and $\max_{i \in \mathcal{V}} \Vert e_{y_j}(\nu) \Vert_{\infty} \leq \frac{1}{2}$ for $\nu \in \{0,1,...,k-1\}$, then we can derive the upper bounds of $\Vert\varsigma(\nu)\Vert_2$. In this way, the upper bounds of $\Vert\Theta(\nu+1)\Vert_2$ can be obtained, which further helps us to derive the non-saturation condition at step $k$.
In other words, if the quantizers are not saturated for all $k \leq k'$, we can obtain $\hat{\varsigma}(1),..., \hat{\varsigma}(k')$ with
\begin{eqnarray}\label{tildeHl}
\hat{\varsigma}(l)=\Vert(\bar{\varsigma}_{1}(l),\bar{\varsigma}_{2}(l),\bar{\varsigma}_{3}(l))^{\text{T}}\Vert_{2}, \quad l\leq k',\notag
\end{eqnarray} 
where $\bar{\varsigma}(l)=\xi^{l}\bar{\varsigma}$ and the elements of the vector 
$\bar{\varsigma}  \in\mathbb{R}^{3}$ is given by: $\bar{\varsigma}_{1}=\frac{1}{2\xi}\eta \pi_{\mathcal{A}}^{\text{T}}\pi_{\mathcal{B}}n\sqrt{m}\beta C $,
$\bar{\varsigma}_{2}=\frac{1}{2\xi}\eta  \kappa_{1} \kappa_{2}\delta_{\mathcal{A}2} n \sqrt{m} \beta C  + \frac{\alpha}{2}\sqrt{mn}\delta_{\mathcal{A}2} \kappa_4 C $,
$\bar{\varsigma}_{3}=
\frac{1}{2\xi}\delta_{\mathcal{B}2}\kappa_3 n\sqrt{m}\beta C(1+ \xi + \eta L )
+ \frac{1}{2}\alpha\delta_{\mathcal{B}2}\kappa_3\kappa_4L\sqrt{mn}C$. 
Note that each element of the vector $\bar{\varsigma}$ is a finite constant.  
We further define the constant $\tilde{\varsigma}$ by: 
\begin{eqnarray}\label{tildevarsi}
\tilde{\varsigma} \triangleq \Vert(\bar{\varsigma}_{1},\bar{\varsigma}_{2},\bar{\varsigma}_{3})^{\text{T}}\Vert_{2}. 
\end{eqnarray} 
Then, we obtain that
$\|\Theta(\iota)\|_2\leq
\|G\|^{\iota}_2 \Vert \Theta(0) \Vert_{2} + \tilde{\varsigma} \sum_{l=0}^{\iota-1}\|G\|^{\iota-1-l}_2\xi^{l}
\leq
\tau \Vert \Theta(0) \Vert_{2} \hat{\rho}^{\iota} + \tilde{\varsigma} \tau \sum_{l=0}^{\iota-2} \hat{\rho}^{\iota-1-l} \xi^{l} + \tilde{\varsigma}\xi^{\iota-1}$
holds,
for $\iota \in \{2,3...,k'+1\}$.

Considering the case $k=k'+1$ ($k' \geq 2$). From (\ref{IMPxmihatx}), we have 
\begin{eqnarray}\label{}
&&\frac{\Vert x_{i}(k)-\hat{x}_{i}(k-1)\Vert_{\infty}}{h(k)}\notag\\
&&\leq
\frac{\sqrt{3}\varphi_{1}}{C \xi^{k}}\left( \tau \Vert \Theta(0) \Vert_{2} \hat{\rho}^{k-1} + \tilde{\varsigma}\tau \sum_{l=0}^{k-3} \hat{\rho}^{k-2-l} \xi^{l} +\tilde{\varsigma}\xi^{k-2} \right) \notag\\
&& ~~~+ \frac{2\alpha n+1}{2\xi} + \frac{ n\eta\beta }{2  \xi^2}
\leq K_x(k) + \frac{1}{2}.
\end{eqnarray}
Similarly, with some tedious calculations, it can also be concluded that $\frac{\Vert y_{i}(k)-\hat{y}_{i}(k-1)\Vert_{\infty}}{h(k-1)}\leq K_y(k) + \frac{1}{2}$.
In summary,
the quantizers will never saturate under the rule (\ref{QuanLevelConAlter}). 
Recalling 
$\Upsilon_{1}(k)$ and $\Upsilon_{2}(k)$ in (\ref{QuanLevelConAlter}), it can be verified that they both can be upper bounded by 
\begin{eqnarray} 
\bar{\Upsilon}= 1+\frac{\tilde{\varsigma}\hat{\rho}}{\xi(\xi-\hat{\rho})\Vert \Theta(0) \Vert_{2}}+\frac{\tilde{\varsigma}}{\xi\tau \Vert \Theta(0) \Vert_{2}}. 
\end{eqnarray} 
Note that $\Vert \Theta(0) \Vert_{2}$, $\hat{\rho}$ and $\tau$ are all some positive constants, $\tilde{\varsigma}$ is a positive constant given in (\ref{tildevarsi}), and $\xi$ 
	 is a constant chosen in the interval $(\hat{\rho},1)$.  
	Hence, $\bar{\Upsilon}$ is a constant, and 
	(\ref{QuanLevelConAlterQQQ}) suffices for the update rule (\ref{QuanLevelConAlter}), 
	 which completes the proof. 
%
\hfill $\blacksquare$



\bibliographystyle{IEEEtran}


\balance
\end{document}